\documentclass[a4,semcolor,psfig,english,11pt,epsf,portrait]{article}
\usepackage{amsmath,amsfonts,latexsym,amscd,amssymb,theorem}
\usepackage{moreverb,rotating,graphics,amsfonts,latexsym,amscd}
\usepackage[ansinew]{inputenc}
\usepackage{epsfig}
\usepackage{amsmath}
\usepackage{amssymb}
\input epsf
\usepackage{psfrag}
\def\R{\hbox{\bf R}}
\def\Z{\hbox{\bf Z}}

\def\o{\overline}

\def\N{\hbox{\bf N}}

\def\g{\gamma}
\def\b{\beta}

\def\a{\alpha}

\def\k{\kappa}

\def\eps{\varepsilon}

\def\l{{\lambda}}

\def\th{\theta}
\def\r{\rho}
\def\t{\tau}

\def\<{\langle}
\def\>{\rangle}
\newcommand{\ba}{\begin{eqnarray}}
\newcommand{\ea}{\end{eqnarray}}

\newtheorem{thm}{Theorem}[section]

\newtheorem{lem}[thm]{Lemma}

\newtheorem{theorem}[thm]{Theorem}
\newtheorem{definition}[thm]{Definition}

\newtheorem{proposition}[thm]{Proposition}

\newtheorem{rem}[thm]{Remark}

\numberwithin{equation}{section}

\renewcommand{\R}{{\mathbb R}}
\renewcommand{\Z}{{\mathbb Z}}
\renewcommand{\N}{{\mathbb N}}

\renewcommand{\g}{\gamma}

\begin{document}

\title{\bf Dynamics of dislocation densities\\ in a bounded channel. Part
  I: smooth solutions to a singular coupled parabolic system.
}
\author{
\normalsize\textsc{ H. Ibrahim $^{*}$, M. Jazar $^{1}$, R. Monneau
\footnote{Cermics, Paris-Est/ENPC, 6 et 8 avenue Blaise Pascal,
Cit\'e Descartes Champs-sur-Marne, 77455 Marne-la-Vall\'ee Cedex 2,
France. E-mails: ibrahim@cermics.enpc.fr,
monneau@cermics.enpc.fr\newline \indent ${}^1$ \hskip-.1cm M. Jazar,
Lebanese University, LaMA-Liban, P.O. Box 826, Tripoli Liban,
mjazar@ul.edu.lb \newline M. Jazar is supported by a grant from
Lebanese University} }} \vspace{20pt}

\maketitle

%%%%%%%%%%%%%%%%%%%%%%%%%%%%%%%%%%%%%%%%%%%%%%%%%%%%%%%%%%%%%%%%%%%%%%%%%%
%%%%%%%%%%%%%%%%%%%%%%%%%%%%%%%%%%%%%%%%%%%%%%%%%%%%%%%%%%%%%%%%%%%%%%%%%%
{\centerline{\small{\bf{Abstract}}} \noindent{\footnotesize{We
study a coupled system of two parabolic equations in one space
dimension. This system is singular because of the presence of one
term with the inverse of the gradient of the solution. Our system
describes an approximate model of the dynamics of dislocation
densities in a bounded channel submitted to an exterior applied
stress. The system of equations is written on a bounded interval
with Dirichlet conditions and requires a special attention to the
boundary. The proof of existence and uniqueness is done under the
use of two main tools: a certain comparison principle on the
gradient of the solution, and a parabolic Kozono-Taniuchi
inequality.}}}

\hfill\break
 \noindent{\small{\bf{AMS Classification: }}} {\small{35K50, 35K40,
     35K55, 42B35,
     42B99.}}\hfill\break
  \noindent{\small{\bf{Key words: }}} {\small{Boundary value problems, parabolic
        systems, nonlinear PDE, $BMO$ spaces,
        logarithmic Sobolev inequality, parabolic Kozono-Taniuchi inequality.}}\hfill\break

%%%%%%%%%%%%%%%%%%%%%%%%%%%%%%%%%%%%%%%%%%%%%%%%%%%%%%%%%%%%%%%
%                                                             %
%                       Section 1                             %
%                                                             %
%                      Introduction                           %
%                                                             %
%%%%%%%%%%%%%%%%%%%%%%%%%%%%%%%%%%%%%%%%%%%%%%%%%%%%%%%%%%%%%%%

\begin{section}{Introduction}\label{sec1}
\subsection{Setting of the problem}
In this paper, we are concerned in the study of the following singular
parabolic system:
\begin{equation}\label{pre_app_model}
\left\{
\begin{aligned}
&\k_{t}=\eps\k_{xx}+\frac{\r_{x}\r_{xx}}{\k_{x}}-\tau\r_{x}\quad&\mbox{on}&\quad
I\times(0,\infty)\\
&\r_{t}=(1+\eps)\r_{xx}-\tau\k_{x}\quad&\mbox{on}&\quad
I\times(0,\infty),
\end{aligned}
\right.
\end{equation}
with the initial conditions:
\begin{equation}\label{ic}
\k(x,0)=\k^{0}(x)\quad\mbox{and}\quad\r(x,0)=\r^{0}(x),
\end{equation}
and the boundary conditions:
\begin{equation}\label{bc}
\left\{
\begin{aligned}
&\k(0,.)=\k^{0}(0)\quad\mbox{and}\quad\k(1,.)=\k^{0}(1),\\
&\r(0,.)=\r(1,.)=0,
\end{aligned}
\right.
\end{equation}
where
$$\eps>0,\quad \tau\in \R,$$
are fixed constants, and
$$I:=(0,1)$$
is the open and bounded interval of $\R$.

The goal is to show the long-time existence and uniqueness of a
smooth solution of (\ref{pre_app_model}), (\ref{ic}) and
(\ref{bc}). Our motivation comes from a problem of studying the
dynamics of dislocation densities in a constrained channel
submitted to an exterior applied stress. In fact, system
(\ref{pre_app_model}) can be seen as an approximate model of the
one described in \cite{GCZ}. This approximate model (presented in
\cite{GCZ} for $\eps=0$) reads:
\begin{equation}\label{sys_theta}
\left\{
\begin{aligned}
&\th^{+}_{t}=\eps\theta^{+}_{xx}+\left[\left(\frac{\th^{+}_{x}-\th^{-}_{x}}{\th^{+}+\th^{-}}
-\tau\right)\th^{+}\right]_{x}\quad&\mbox{on}&\quad I\times(0,\infty),\\
&\th^{-}_{t}=\eps\theta^{-}_{xx}-\left[\left(\frac{\th^{+}_{x}-\th^{-}_{x}}{\th^{+}+\th^{-}}
-\tau\right)\th^{-}\right]_{x}\quad&\mbox{on}&\quad
I\times(0,\infty),
\end{aligned}
\right.
\end{equation}
with $\t$ representing the exterior stress field. System
(\ref{sys_theta}) can be deduced from (\ref{pre_app_model}), by
spatially differentiating (\ref{pre_app_model}), and by considering
\begin{equation}\label{part1:thekr}
\r^{\pm}_{x}=\theta^{\pm},\quad \r=\r^{+}-\r^{-},\quad
\k=\r^{+}+\r^{-},
\end{equation}
which explains the presence of the factor $(1+\eps)$ in the second
equation of (\ref{pre_app_model}). Here $\th^{+}$ and $\theta^{-}$
represent the densities of the positive and negative dislocations
respectively (see \cite{nab, HL} for a physical study of
dislocations).

The part II of this work will be presented in \cite{IJM_PII}. There,
we will show some kind of convergence of the solution
$(\r^{\eps},\k^{\eps})$ as $\eps \rightarrow 0$.

\subsection{Statement of the main result}
The main result of this paper is:
\begin{theorem}\textit{\textbf{(Existence and uniqueness of a
      solution).}}\label{mainresult0}
Let $\r^{0}$, $\k^{0}$ satisfying:
\begin{equation}\label{MT1}
\r^{0},\k^{0}\in
C^{\infty}(\bar{I}),\quad \r^{0}(0)=\r^{0}(1)=\k^{0}(0)=0,\quad \k^{0}(1)=1,
\end{equation}
\begin{equation}\label{misi2z}
\left\{
\begin{aligned}
&(1+\eps)\r^{0}_{xx}=\t\k^{0}_{x}\quad &\mbox{on}&\quad \partial I\\
&(1+\eps)\k^{0}_{xx}=\t\r^{0}_{x}\quad &\mbox{on}&\quad \partial I,
\end{aligned}
\right.
\end{equation}
and
\begin{equation}\label{MT3}
\k^{0}_{x}>|\r^{0}_{x}|\quad \mbox{on}\quad \bar{I}.
\end{equation}
Then there exists a unique global solution $(\r,\k)$ of system
(\ref{pre_app_model}), (\ref{ic}) and (\ref{bc}) satisfying
\begin{equation}\label{misi3z}
(\r,\k)\in C^{3+\a,\frac{3+\a}{2}}(\bar{I}\times [0,\infty))\cap
C^{\infty}(\bar{I}\times (0,\infty)),\quad \forall \a \in (0,1).
\end{equation}
Moreover, this solution also satisfies :
\begin{equation}\label{MT5}
\k_{x}>|\r_{x}| \quad \mbox{on} \quad  \bar{I}\times [0,\infty).
\end{equation}
\end{theorem}
\begin{rem}
Conditions  (\ref{misi2z}) are natural here. Indeed, the regularity
(\ref{misi3z}) of the solution of  (\ref{pre_app_model}) with the
boundary conditions (\ref{ic}) and (\ref{bc}) imply in particular
(\ref{misi2z}).
\end{rem}
\begin{rem}
Remark that the choice $\k^{0}(0)=0$ and $\k^{0}(1)=1$ does not
reduce the generality of the problem, because equation
(\ref{pre_app_model}) does not see the constants and has the
following invariance: if $(\r,\k)$ is a solution, then $(\lambda \r,
\lambda \kappa)$ is also a solution for any $\lambda \in \R$.
\end{rem}
\subsection{Brief review of the literature}
To our knowledge, systems of equations involving the singularity in
$1/\k_{x}$ as in (\ref{pre_app_model}) has not been directly handled
elsewhere in the literature. However, parabolic problems involving
singular terms have been widely studied in various aspects. Fast
diffusion equations:
$$
u_{t}-\Delta u^{m}=0,\quad 0<m<1,
$$
are examined, for instance, in \cite{vazquez_carrillo03,
vazquez_cha02, DiBook}. These equations are singular at points where
$u=0$. In dimension 1, setting $u=v_{x}$ we get, up to a constant of
integration:
$$v_{t}-m v_{x}^{m-1}v_{xx}=0$$
which makes appear a singularity like $1/v_{x}$. Other class of
singular parabolic equations are for instance of the form:
\begin{equation}\label{lit3}
u_{t}=u_{xx}+\frac{b}{x}u_{x},
\end{equation}
where $b$ is a certain constant. Such an equation is related to
axially symmetric problems and also occurs in probability theory
(see \cite{ChanKaper, Antonio}). An important type of equations that
can be indirectly related to our system are semilinear parabolic
equations:
\begin{equation}\label{sale2HP}
u_{t}=\Delta u + |u|^{p-1}u,\quad p>1.
\end{equation}
Many authors have studied the blow-up phenomena for solutions of the above
equation (see for instance \cite{MerleZaag, SouGuo}).
Equation (\ref{sale2HP}) can be somehow related to the
first equation of (\ref{pre_app_model}), but with a singularity of the
form ${1}/{\k}$. This can be formally seen if we first suppose that
$u\geq 0$, and then we apply the following
change of variables $u=1/v$. In this case, equation (\ref{sale2HP})
becomes:
$$
v_{t}=\Delta v -\frac{2|\nabla v|^{2}}{v}-v^{2-p},
$$
and hence if $p=3$, we obtain:
\begin{equation}\label{sale2HP2}
v_{t}=\Delta v -\frac{1}{v}(1+2|\nabla v|^{2}).
\end{equation}
Since the solution $u$ of (\ref{sale2HP}) may blow-up at a finite
time $t=T$, then $v$ may vanishes at $t=T$, and therefore equation
(\ref{sale2HP2}) faces similar singularity to that of the first
equation of (\ref{pre_app_model}), but in terms of the solution
$v$ instead of $v_x$.

\subsection{Strategy of the proof}
The existence and uniqueness is made by using a fixed point argument
after a slight artificial modification in the denominator $\k_{x}$ of
the first equation of (\ref{pre_app_model}) in order to avoid dividing
by zero.  We will first show the short
time existence, proving in particular that
$$
\k_{x}(x,t)\geq \sqrt{\g^{2}(t)+\r_{x}^{2}(x,t)}\geq \g(t)>0,
$$
for some well chosen initial data and a suitable function $\g(t)>0$.
The only, but dangerous, inconvenience is that the function $\g$
depends strongly on $\|\r_{xxx}(.,t)\|_{L^{\infty}(I)}$, roughly
speaking:
\begin{equation}\label{cc:et}
\g^{'}\simeq -\|\r_{xxx}\|_{L^{\infty}(I)} \g.
\end{equation}
Let us mention that one of the key points here is that
$\left|\frac{\r_{x}}{\k_{x}}\right|\leq 1$ which somehow linearize
the first equation of (\ref{pre_app_model}). Nevertheless, standard
Sobolev and Hölder estimates for the parabolic system
(\ref{pre_app_model}) are not good enough to bound
$\|\r_{xxx}\|_{L^{\infty}(I)}$ in order to prevent $\g$ (and as a
consequence $\k_{x}$) from vanishing. On the contrary, a Sobolev
logarithmic estimate (see Section \ref{sec2}, the parabolic
Kozono-Taniuchi inequality, Theorem \ref{KTRineq}) can be used in
order to obtain a sharp bound of $\|\r_{xxx}\|_{L^{\infty}(I)}$ of
the form
$$\|\r_{xxx}\|_{L^{\infty}(I)}\leq E\left(1+\log^{+}\frac{1}{\g}\right),$$
where $E$ is an exponential function in time. This allows, with
(\ref{cc:et}), to show that the function $\g>0$ does not vanish in
finite time. After that, due to some \textit{a priori} estimates, we
can prove the global time existence.

\subsection{Organization of the paper}
This paper is organized as follows: in Section \ref{sec2}, we
present the tools needed throughout this work, this includes a brief
recall on the $L^{p}$, $C^{\a}$ and  the $BMO$ theory for parabolic
equations. In Section \ref{sec3}, we show a comparison principle
associated to (\ref{pre_app_model}) that will play a crucial rule in
the long time existence of the solution as well as the positivity of
$\k_{x}$. In Section \ref{sec4}, we present a result of short time
existence, uniqueness and regularity of a solution $(\r,\k)$ of
(\ref{pre_app_model}). Section \ref{sec5} is devoted to give some
exponential bounds on $(\r,\k)$. In Section \ref{sec5.1}, we show a
control of the $W^{2,1}_{2}$ norm of $\r_{xxx}$. In a similar way,
we show a control of the $BMO$ norm of $\r_{xxx}$ in Section
\ref{sec5.2}. In Section \ref{sec5.3}, we use a parabolic
Kozono-Taniuchi inequality to control the $L^{\infty}$ norm of
$\r_{xxx}$. In Section \ref{sec6}, we prove our main result: Theorem
\ref{mainresult0}. Finally, Sections \ref{secA} and \ref{secB} are
appendices where we present the proofs of some technical results.
\end{section}

%%%%%%%%%%%%%%%%%%%%%%%%%%%%%%%%%%%%%%%%%%%%%%%%%%%%%%%%%%%%%%
%                                                            %
%                                                            %
%                         Section 2                          %
%                                                            %
%                           Tools                            %
%                                                            %
%%%%%%%%%%%%%%%%%%%%%%%%%%%%%%%%%%%%%%%%%%%%%%%%%%%%%%%%%%%%%%

\begin{section}{Tools: theory of parabolic equations}\label{sec2}

We start with some basic notations and terminology:\\

\noindent {\bf {\underline {Abridged notation.}}}\\

\noindent $\bullet$ $I_{T}$ is the cylinder $I\times (0,T)$; $\bar{I}$
is the closure of $I$;  $\overline{I_{T}}$
is the closure of $I_{T}$; $\partial I$ is the boundary of $I$.\\

\noindent $\bullet$ $\|.\|_{L^{p}(\Omega)}=\|.\|_{p,\Omega}$,\;
$\Omega$
is an open set,\; $p\geq 1$.\\

\noindent $\bullet$ $S_{T}$ is the lateral boundary of $I_{T}$, or
more
precisely, $S_{T}=\partial I\times (0,T)$.\\

\noindent $\bullet$ $\partial^{p}I_{T}$ is the parabolic boundary of
$I_{T}$, i.e. $\partial^{p}I_{T}= \overline{S_{T}} \cup (I\times
\{t=0\})$.\\

\noindent $\bullet$ $D^{s}_{y}u=\frac{\partial^{s}u}{\partial y^{s}}$,\; $u$ is a function
depending on the parameter $y$,\;  $s\in \N.$\\

\noindent $\bullet$ $[l]$ is the floor part of $l\in \R$.\\

\noindent $\bullet$ $Q_{r}=Q_{r}(x_{0},t_{0})$ is the lower parabolic
cylinder
given by: $$Q_{r}=\{(x,t);\; |x-x_{0}|<r,\;
t_{0}-r^{2}<t<t_{0}\},\;r>0, \;(x_{0},t_{0})\in I_{T}.$$

\noindent $\bullet$ $|\Omega|$ is the $n$-dimensional Lebesgue
measure of the open set $\Omega\subset \R^{n}$.\\

\noindent $\bullet$ $m_{\Omega}(u)=\frac{1}{|\Omega|}\int_{\Omega} u$ is
the average integral of the $u\in L^{1}(\Omega)$ over $\Omega\subset
\R^{n}$.

\subsection{$L^{p}$  and $C^{\a}$ theory of parabolic equations}

A major part of this work deals with the following typical problem in
parabolic theory:
\begin{equation}\label{typpara}
\left\{
\begin{aligned}
&u_{t}=\eps u_{xx}+f\quad&\mbox{on}&\quad I_T\\
&u(x,0)=\phi \quad&\mbox{on}&\quad I\\
&u=\Phi\quad&\mbox{on}&\quad \partial I\times (0,T),
\end{aligned}
\right.
\end{equation}
where $T>0$ and $\eps>0$. A wide literature on the existence and
uniqueness of solutions of (\ref{typpara}) in different function
spaces could be found for instance in \cite{LSU}, \cite{afriedman}
and \cite{lieberman}. We will deal mainly with two types of
spaces:\\

\noindent {\bf The Sobolev space} $W^{2,1}_{p}(I_{T})$, $1<p<\infty$
which
  is the Banach space consisting of the elements in
$L^{p}(I_{T})$ having generalized derivatives of the form
$D^{r}_{t}D^{s}_{x}u$, with $r$ and $s$ two non-negative integers
satisfying the inequality $2r+s\leq 2$, also in $L^{p}(I_{T})$. The
norm in this space is defined as
$\|u\|_{W^{2,1}_{p}(I_{T})}=\sum_{i=0}^{2}\sum_{2r+s=i}\|D^{r}_{t}D^{s}_{x}u\|_{p,I_{T}}$.\\

\noindent {\bf The Hölder spaces} $C^{\ell}(\bar I)$  and
  $C^{\ell,\ell/2}(\overline{I_{T}})$, $\ell>0$ a nonintegral positive
  number. We do not recall the definition of the space $C^{\ell}(\bar
  I)$ which is very standard. The Hölder space $C^{\ell,\ell/2}(\overline{I_{T}})$ is the
  Banach space of functions $v(x,t)$ that are continuous in
  $\overline{I_{T}}$, together with all derivatives of the form
  $D^{r}_{t}D^{s}_{x}v$ for $2r+s<\ell$, and have a finite norm
$|v|_{I_{T}}^{(\ell)}=\<v\>^{(\ell)}_{I_{T}}+\sum_{j=0}^{[\ell]}\<v\>_{I_{T}}^{(j)}$,
where
$$\<v\>_{I_{T}}^{(0)}=|v|^{(0)}_{I_{T}}=\|v\|_{\infty,I_{T}},\quad
\<v\>_{I_{T}}^{(j)}=\sum_{2r+s=j}|D^{r}_{t}D^{s}_{x}v|^{(0)}_{I_{T}},\quad
\<v\>^{(\ell)}_{I_{T}}=\<v\>_{x,I_{T}}^{(\ell)}+\<v\>_{t,I_{T}}^{(\ell/2)},$$
and
$$\<v\>_{x,I_{T}}^{(\ell)}=\sum_{2r+s=[\ell]}\<D^{r}_{t}D^{s}_{x}v\>_{x,I_{T}}^{(\ell-[\ell])},\quad
\<v\>_{t,I_{T}}^{(\ell/2)}=\sum_{0<\ell-2r-s<2}\<D^{r}_{t}D^{s}_{x}v\>_{t,I_{T}}^{\left(\frac{\ell-2r-s}{2}\right)},$$
with
$$\<v\>_{x,I_{T}}^{(\a)}=\inf\{c;\; |v(x,t)-v(x{'},t)|\leq
c|x-x{'}|^{\a},\;(x,t), (x',t)\in \overline{I_{T}}\},\quad 0<\a<1,
$$
$$
\<v\>_{t,I_{T}}^{(\a)}=\inf\{c;\; |v(x,t)-v(x,t{'})|\leq
c|t-t{'}|^{\a},\;(x,t), (x,t')\in \overline{I_{T}}\},\quad 0<\a<1.
$$
The above definitions could be found in details in \cite[Section 1]{LSU}.
Now, we write down the compatibility conditions of order $0$ and
$1$. These compatibility conditions concern the given data $\phi$,
$\Phi$ and $f$
of problem (\ref{typpara}).\\

\noindent {\bf Compatibility condition of order 0.} Let $\phi\in
C(\bar{I})$ and $\Phi\in C(\overline{S_{T}})$. We say that the compatibility
condition of order $0$ is satisfied if
\begin{equation}\label{tibi0}
\phi\big|_{\partial I}=\Phi\big|_{t=0}.
\end{equation}

\noindent {\bf Compatibility condition of order 1.} Let $\phi\in
C^{2}(\bar{I})$, $\Phi\in C^{1}(\overline{S_{T}})$ and $f\in C(\overline{I_{T}})$. We
say that the compatibility condition of order $1$ is satisfied if
(\ref{tibi0}) is satisfied and in addition we have:
\begin{equation}\label{tibi1}
(\eps\phi_{xx}+f)\big|_{\partial I}=\frac{\partial \Phi}{\partial t}\Big|_{t=0}.
\end{equation}
We state two results of existence and uniqueness adapted to our
special problem. We begin by presenting the solvability of parabolic
equations in Hölder spaces.
\begin{theorem}\textit{\textbf{(Solvability in Hölder spaces, \cite[Theorem
5.2]{LSU}).}}\label{holder_sol} Suppose \break $0<\alpha<2$, a
non-integral number. Then for any
$$\phi\in C^{2+\alpha}(\bar{I}), \quad \Phi\in C^{1+\alpha/2}(\overline{S_{T}})\quad \mbox{and}\quad f\in
C^{\alpha,\alpha/2}(\overline{I_{T}})$$ satisfying the compatibility
condition of order $1$ (see (\ref{tibi0}) and (\ref{tibi1})),
problem (\ref{typpara}) has a unique solution $u\in
C^{2+\a,1+\a/2}(\overline{I_{T}})$ satisfying the following
inequality:
\begin{equation}\label{cor:hol_ineq}
|u|_{I_{T}}^{(2+\a)}\leq c^{H}
\left(|f|_{I_{T}}^{(\a)}+|\phi|_{I}^{(2+\a)}+|\Phi|_{{S}_{T}}^{(1+\a/2)}\right),
\end{equation}
for some $c^{H}=c^{H}(\eps, \a, T)>0$.
\end{theorem}
\begin{rem}\label{notwritten}\textit{\textbf{(Estimating the term $c^{H}$ of (\ref{cor:hol_ineq})).}}
The constant appearing in the above Hölder estimate
(\ref{cor:hol_ineq}) can be estimated, using some time iteration, as
$ c^{H}(\eps, \a, T)\leq e^{c(T+1)}$, where $c=c(\eps,\a)>0$ is a
positive constant.
\end{rem}
We now present the solvability in Sobolev spaces. Recall the norm
of fractional Sobolev spaces. If $f\in W^{s}_{p}(a,b)$, $s>0$ and
$1<p<\infty$, then
\begin{equation}\label{frac_norm}
\|f\|_{W^{s}_{p}(a,b)}=\|f\|_{W^{[s]}_{p}(a,b)}+\left(\int_{a}^{b}\int_{a}^{b}
\frac{|f^{([s])}(x)-f^{([s])}(y)|^{p}}{|x-y|^{1+(s-[s])p}}\right)^{1/p}.
\end{equation}

\begin{theorem}\textit{\textbf{(Solvability in Sobolev spaces, \cite[Theorem
    9.1]{LSU}).}}\label{classpara}
Let $p>1$,\break $\eps>0$ and $T>0$. For any $f\in L^{p}(I_{T})$,
$\phi\in W^{2-2/p}_{p}(I)$ and $\Phi\in W^{1-1/2p}_{p}(S_{T})$, with
\break $p\neq 3/2$ ($p=3/2$ is called the \textbf{singular} index)
satisfying in the case $p>3/2$ the compatibility condition of order
zero (see (\ref{tibi0})), there exists a unique solution $u\in
W^{2,1}_{p}(I_{T})$ of (\ref{typpara}) satisfying the following
estimate:
\begin{equation}\label{parabolic_est}
\|u\|_{W^{2,1}_{p}(I_{T})}\leq
c\,\left(\|f\|_{p,I_{T}}+\|\phi\|_{W^{2-2/p}_{p}(I)}+\|\Phi\|_{W^{1-1/2p}_{p}(S_{T})}\right),
\end{equation}
for some $c=c(\eps, p, T)>0$.
\end{theorem}
For a better understanding of the spaces stated in the above two
theorems, especially fractional Sobolev spaces, we send the reader
to \cite{adams} or \cite{LSU}. The dependence of the constant $c$
of Theorem \ref{classpara} on the variable $T$ will be of notable
importance and this what is emphasized by the next lemma.
\begin{lem}\textit{\textbf{(The constant $c$ given by (\ref{parabolic_est}): case
    $\phi=0$ and $\Phi=0$).}}
\label{prec_estpara} Under the same hypothesis of Theorem
\ref{classpara}, with $\phi=0$ and $\Phi=0$, the estimate
(\ref{parabolic_est}) can be written as:
\begin{equation}\label{parabolicinq}
\frac{\|u\|_{p,I_{T}}}{T}+\frac{\|u_{x}\|_{p,I_{T}}}{\sqrt{T}}+
\|u_{xx}\|_{p,I_{T}}+\|u_{t}\|_{p,I_{T}}\leq
{c}\|f\|_{p,I_{T}},
\end{equation}
where $c=c(\eps, p)>0$ is a positive constant depending only on $p$ and $\eps$.
\end{lem}
The proof of this lemma will be done in  Appendix A. Moreover, We will frequently
make use of the following two lemmas also depicted  from \cite{LSU}.
\begin{lem}\textit{\textbf{(Sobolev embedding in Hölder spaces, \cite[Lemma
    3.3]{LSU}).}}\label{lemmalsu1}\\
(i) (Case $p>3$). For any function $u \in W^{2,1}_{p}(I_{T})$, if
$\alpha=1-{3}/{p}>0$, i.e. $p>3$, then $u\in
C^{1+\a,\frac{1+\a}{2}}(\overline{I_{T}})$ with
$|u|^{(1+\a)}_{I_{T}}\leq c\|u\|_{W^{2,1}_{p}(I_{T})}$,
$c=c(p,T)>0$. However, in terms of $u_{x}$, we have that $ u_{x}\in
C^{\alpha,\alpha/2}(\overline{I_{T}})$ satisfies the following
estimates:
\begin{equation}\label{cor:estH1}
\|u_{x}\|_{\infty,{I}_{T}}\leq c\left\{
    \delta^{\alpha}(\|u_{t}\|_{p,I_{T}}+\|u_{xx}\|_{p,I_{T}})+\delta^{\alpha-2}\|u\|_{p,I_{T}}\right\},\quad c=c(p)>0,
\end{equation}
$$
\<u_{x}\>^{(\alpha)}_{I_{T}}\leq c\left\{
    \|u_{t}\|_{p,I_{T}}+\|u_{xx}\|_{p,I_{T}}+\frac{1}{\delta^{2}}\|u\|_{p,I_{T}}\right\},\quad c=c(p)>0.
$$
(ii) (Case $p>3/2$). If $u\in W^{2,1}_{p}(I_{T})$ with $p>3/2$, then
$u\in C(\overline{I_{T}})$, and we have the following estimate:
\begin{equation}\label{almz}
\|u\|_{\infty,{I}_{T}}\leq c\left\{
    \delta^{2-3/p}(\|u_{t}\|_{p,I_{T}}+\|u_{xx}\|_{p,I_{T}})+\delta^{-3/p}\|u\|_{p,I_{T}}\right\},\quad c=c(p)>0.
\end{equation}
In the above two cases $\delta=\min\{1/2,\sqrt{T}\}$.
\end{lem}
\begin{lem}\textit{\textbf{(Trace of functions in $W^{2,1}_{p}(I_{T})$,
      \cite[Lemma 3.4]{LSU}).}}\label{lemmalsu}
If $u\in W^{2,1}_{p}(I_{T})$, $p>1$, then for $2r+s<2-2/p$, we have
$D^{r}_{t}D^{s}_{x}u{\big |_{t=0}}\in W^{2-2r-s-2/p}_{p}(I)$ with
$$
\|u\|_{W^{2-2r-s-2/p}_{p}(I)}\leq
c(T)\|u\|_{W^{2,1}_{p}(I_{T})}.
$$
In addition, for $2r+s<2-1/p$, we have
$D^{r}_{t}D^{s}_{x}u{\big|_{\overline{S_{T}}}}\in
W^{1-r-s/2-1/2p}_{p}(\overline{S_{T}}) $ with
$$
\|u\|_{W^{1-r-s/2-1/2p}_{p}\left(\overline{S_{T}}\right)}\leq
c(T)\|u\|_{W^{2,1}_{p}(I_{T})}.
$$
\end{lem}
A useful technical lemma will now be presented. The proof of this lemma
will be done in Appendix A.
\begin{lem}\textit{\textbf{($L^{\infty}$ control of the spatial
      derivative).}}\label{cor:con_hol_sob}
Let $p>3$ and let $0<T\leq 1/4$ (this condition is taken for
simplification). Then for every $u\in W^{2,1}_{p}(I_{T})$ with $u=0$
on $\partial^{p}(I_{T})$ in the trace sense (see Lemma
(\ref{lemmalsu})), there exists a constant $c(T,p)>0$ such that
$$
\|u_{x}\|_{\infty,I_{T}}\leq c(T,p)\|u\|_{W^{2,1}_{p}(I_{T})},\quad
\mbox{with}\quad c(T,p)=c(p)T^{\frac{p-3}{2p}}\rightarrow
0\,\,\mbox{as}\,\,T\rightarrow 0.
$$
\end{lem}
\subsection{$BMO$ theory for parabolic equation}
A very useful tool in this paper is the limit case of the $L^{p}$
theory, $1<p<\infty$, for parabolic equations, which is the $BMO$
theory. Roughly speaking, if the function $f$ appearing in
(\ref{typpara}) is in the $L^{p}$ space for some $1<p<\infty$, then
we expect our solution $u$ to have $u_{t}$ and $u_{xx}$ also in
$L^{p}$. This is no longer valid in the limit case, i.e. when
$p=\infty$. In this case, it is shown that the solution $u$ of the
parabolic equation have $u_{t}$ and $u_{xx}$ in the
parabolic/anisotropic $BMO$ (bounded mean oscillation) space  that
is convenient to present some of its related theories.
\begin{definition}\textit{\textbf{(Parabolic/Anisotropic $BMO$
      spaces).}}\label{def_bmo}
A function $u\in L^{1}_{loc}(I_{T})$ is said to be of bounded mean
oscillation, $u\in BMO(I_{T})$, if the quantity
$$
\sup_{Q_{r}\subset I_{T}}\left(\frac{1}{|Q_{r}|}\int_{Q_r}\left|u-m_{Q_{r}}(u)\right|\right)
$$
is finite. Here the supremum is taken over all parabolic lower
cylinders $Q_{r}$.
\end{definition}
\begin{rem}
The parabolic $BMO(I_{T})$ space, which will be refereed, for
simplicity, as the $BMO(I_{T})$ space, and sometimes, where there is
no confusion, as $BMO$ space, is a Banach space (whose elements are
defined up to an additive constant) equipped with the norm
$$
\|u\|_{BMO(I_{T})}=\sup_{Q_{r}\subset
  I_{T}}\left(\frac{1}{|Q_{r}|}\int_{Q_r}\left|u-m_{Q_{r}}(u)\right|\right).
$$
\end{rem}
We move now to the two main theorems of this subsection, the $BMO$
theory for parabolic equations, and the Kozono-Taniuchi parabolic
type inequality. To be more precise, we have the following:

\begin{theorem}\textit{\textbf{($BMO$ theory for parabolic
    equations in the periodic case).}}\label{bmo_theory}
Take $0<T_{1}\leq T$. Consider the following Cauchy problem:
\begin{equation}\label{cau:1}
\left\{
\begin{aligned}
& u_{t}=\eps u_{xx}+f\quad &\mbox{on}& \quad \R\times (0,T),\\
& u(x,0)=0.
\end{aligned}
\right.
\end{equation}
If $f\in L^{\infty}(\R\times (0,T))$ and $f$ is a $2I$-periodic
function in space, i.e. $f(x+2,t)=f(x,t)$, then there exists a
unique solution $u\in BMO(\R\times (0,T))$ of (\ref{cau:1})
with\break  $u_{t}$, $u_{xx}\in BMO(\R\times (0,T))$. Moreover,
there exists $c>0$ that may depend on $T_{1}$ but independent of $T$
such that:
\begin{equation}\label{cau:2}
\|u_{t}\|_{BMO(\R\times (0,T))}+\|u_{xx}\|_{BMO(\R\times (0,T))}\leq
c \big[\|f\|_{BMO(\R\times (0,T))}+m_{2I\times (0,T)}(|f|)\big].
\end{equation}
\end{theorem}
The proof of this theorem will be presented in Appendix B. Our next
tool (see Theorem \ref{KTRineq}) shows an estimate involving
parabolic $BMO$ spaces. This estimate is a control of the
$L^{\infty}$ norm of a given function by its $BMO$ norm and the
logarithm of its norm in a certain Sobolev space. It can also be
considered as the parabolic version on a bounded domain $I_{T}$ of
the Kozono-Taniuchi inequality (see \cite{KT}) that we recall here.
\begin{theorem}\textit{\textbf{(The Kozono-Taniuchi inequality in the elliptic
    case, \cite[Theorem 1]{KT}).}}\label{koz-tan}
Let $1<p<\infty$ and let $s>n/p$. There is a constant $C=C(n,p,s)$
such that, for all $f\in W_{p}^{s}(\R^{n})$, the following estimate
holds:
\begin{equation}\label{koztanest}
\|f\|_{\infty,\R^{n}}\leq
C\left(1+\|f\|_{BMO_e(\R^{n})}\left(1+\log^{+}\|f\|_{W_{p}^{s}(\R^{n})}\right)
\right),\,\,\log^{+}=\max(0,\log).
\end{equation}
\end{theorem}
\begin{rem}
It is worth mentioning that the $BMO_e$ norm appearing in
(\ref{koztanest}) is the elliptic $BMO_e$ norm, i.e. the one where
the supremum is taken over ordinary balls
$$B_{r}(X_0)=\{X\in\R^{n};\; |X-X_{0}|<r \}.$$
\end{rem}
The original type of the logarithmic Sobolev inequality was found in
\cite{BreGal, BreWai} (see also \cite{Engler}), where the authors investigated
the relation between $L^{\infty}$, $W^{k}_{r}$ and $W^{s}_{p}$ and
proved that there holds the embedding
$$\|u\|_{L^{\infty}(\R^{n})}\leq C
\left(1+\log^{\frac{r-1}{r}}\left(1+\|u\|_{W^{s}_{p}(\R^{n})}\right)\right),\quad
sp>n$$ provided $\|u\|_{W^{k}_{r}}\leq 1$ for $kr=n$. This
estimate was applied to prove existence of global solutions to the
nonlinear Schrödinger equation (see \cite{BreGal,HayWol}).

In our work, we need to have an estimate similar to (\ref{koztanest}),
but for the parabolic $BMO$ space and on the bounded domain
$I_{T}$. This will be essential, on one hand, to show a suitable
positive lower bound of $\k_{x}$ ($\k$ given by Theorem \ref{mainresult0}),
and on the other hand, to show the long time existence of our
solution. Indeed, there is a similar inequality and this is
what will be illustrated by the next theorem.
\begin{theorem}\textit{\textbf{(A parabolic Kozono-Taniuchi
      inequality, \cite[Appendix B2]{IJM_PI}, \cite{Ib_Mo08}).}}\label{KTRineq}
Let $v\in W^{2,1}_{2}(I_{T})$, then there exists a constant
$c=c(T)>0$ such that the estimate holds
\begin{equation}\label{kte_est}
\|v\|_{\infty,I_{T}}\leq c\|v\|_{\overline{BMO}(I_{T})}\left(1+\log^{+}\|v\|_{W^{2,1}_{2}(I_{T})}\right),
\end{equation}
where $\overline{BMO}(I_{T})=BMO(I_{T})\cap L^{1}(I_{T})$, and for $v\in
\overline{BMO}(I_{T})$,
$$\|v\|_{\overline{BMO}(I_{T})}=\|v\|_{BMO(I_{T})}+
\|v\|_{L^{1}(I_{T})}.$$ 
\end{theorem}
This inequality is first shown over $\R_{x}\times \R_{t}$, then it
is deduced over $I_{T}$.
\end{section}

%%%%%%%%%%%%%%%%%%%%%%%%%%%%%%%%%%%%%%%%%%%%%%%%%%%%%%%%%%%%%%%%%%%
%                                                                 %
%                        Section 3                                %
%                                                                 %
%                   Comparison principle                          %
%                                                                 %
%%%%%%%%%%%%%%%%%%%%%%%%%%%%%%%%%%%%%%%%%%%%%%%%%%%%%%%%%%%%%%%%%%%

\begin{section}{A comparison principle}\label{sec3}
\begin{proposition}\textit{\textbf{(A comparison principle for system
    (\ref{pre_app_model})).}}\label{max_p} Let
$$(\r,\k)\in \big(C^{3+\a,\frac{3+\a}{2}}\left(\overline{I_{T}}\right)\big)^{2}\quad
\mbox{for some}\quad 0<\a<1,$$ be a solution of
(\ref{pre_app_model}), (\ref{ic}) and (\ref{bc}) with $\k_{x}>0$,
and the initial conditions $\r^{0}$, $\k^{0}$ satisfying:
\begin{equation}\label{minko0r0x}
\k_{x}^{0}\geq\sqrt{\g_{0}^{2}+(\r_{x}^{0})^{2}} \quad \mbox{on}\quad
I, \quad \g_{0}\in (0,1).
\end{equation}
Choose $\beta=\beta(\eps,\tau)>0$ large enough. Let the function $\g(t)$
satisfies:
\begin{equation}\label{from_cras}
\left\{
\begin{aligned}
& \frac{\g'(t)}{\g(t)}\leq
-\left(c_{0}+\|\r_{xxx}(.,t)\|_{L^{\infty}(I)}\right),\quad
c_{0}=c_{0}(\eps,\beta,\tau),\\
& \g(0)=\g_{0}/2.
\end{aligned}
\right.
\end{equation}
Define $\overline{M}(x,t):=\cosh(\beta
x)\{\kappa_{x}(x,t)-\sqrt{\g^{2}(t)+(\r_{x}(x,t))^{2}}\}$ for
$(x,t)\in \overline{I_{T}}$ and $T>0$. Then $\displaystyle
\overline{m}(t):=\min_{x\in I}\overline{M}(x,t)$ satisfies
$\overline{m}(t)\geq \g^{2}(t)$ for all $t\in [0,T]$. In particular,
we have
\begin{equation}\label{hra123}
\k_{x}(x,t)\geq \sqrt{\g^{2}(t)+\r^{2}_{x}(x,t)},\quad\mbox{in}\quad
\overline{I_{T}}.
\end{equation}
\end{proposition}
\noindent {\bf Proof.} Throughout the proof, we will extensively use the
following notation:
$$
G_{a}(y)=\sqrt{a^{2}+y^{2}},\quad a,y\in \R.
$$
Without loss of generality (up to a change of variables in $(x,t)$ and a
re-definition of $\t$), assume in the proof that
$$I=(-1,1).$$
Define the quantity $M$ by:
$$
M(x,t)=\k_{x}(x,t)-G_{\g(t)}(\r_{x}(x,t)),\quad (x,t)\in
\overline{I_{T}},
$$
$\g(t)>0$ is a function to be determined.
 The proof could be divided into five steps.\\

\noindent {\bf Step 1. }\textsf{(Partial differential inequality
  satisfied by $M$)}\\

\noindent We first do the following computations on $I_{T}$:
\begin{equation}\label{compu_1}
M_{t}=\k_{xt}-G^{'}_{\g}(\r_{x})\r_{xt}-\frac{\g\g^{'}}{\sqrt{\g^{2}+\r_{x}^{2}}},
\end{equation}
\begin{equation}\label{compu_2}
M_{x}=\k_{xx}-G^{'}_{\g}(\r_{x})\r_{xx},\quad
M_{xx}=\k_{xxx}-G^{''}_{\g}(\r_{x})\r^{2}_{xx}-G^{'}_{\g}(\r_{x})\r_{xxx}.
\end{equation}
Deriving (\ref{pre_app_model}) with respect to $x$, we deduce that
\begin{equation}\label{compu_3and4}
\left\{
\begin{aligned}
&\k_{xt}=\eps\k_{xxx}+\frac{\r_{xx}^{2}}{\k_{x}}+\frac{\r_{x}\r_{xxx}}{\k_{x}}-
\frac{\r_{x}\r_{xx}\k_{xx}}{\k_{x}^{2}}-\tau\r_{xx},\\
&\r_{xt}=(1+\eps)\r_{xxx}-\tau\k_{xx}.
\end{aligned}
\right.
\end{equation}
We set
$$\Gamma=\frac{\g\g^{'}}{\sqrt{\g^{2}+\r_{x}^{2}}},\quad  F_{\g}(y)=y-\g\arctan(y/\g).$$
Doing again some direct computations, and using (\ref{compu_1}),
(\ref{compu_2}) and (\ref{compu_3and4}), we obtain
\begin{equation}\label{M_t}
\begin{aligned}
&M_{t}=\eps M_{xx}+\left( \tau
G^{'}_{\g}(\r_{x})-\frac{\r_{x}\r_{xx}}{\k_{x}^{2}}\right)M_{x}+\left(
\frac{\r_{xx}^{2}}{\k_{x}^{2}}-\frac{\r_{xxx}G^{'}_{\g}(\r_{x})}{\k_{x}}\right)M\\
&+\eps
G^{''}_{\g}(\r_{x})\r_{xx}^{2}+\frac{\r_{xx}^{2}}{\k_{x}^{2}}[G_{\g}(\r_{x})-G^{'}_{\g}(\r_{x})\r_{x}]-
\tau(1-F^{'}_{\g}(\r_{x}))\r_{xx}-\Gamma.
\end{aligned}
\end{equation}
Using Young's inequality $2ab\leq a^{2}+b^{2}$, we have:
\begin{equation}\label{young}
\frac{\tau\g^{2}|\r_{xx}|}{\g^{2}+\r_{x}^{2}}\leq
\frac{\eps\g^{2}\r_{xx}^{2}}{(\g^{2}+\r_{x}^{2})^{3/2}}+\frac{\g^{2}\tau^{2}}{4\eps\sqrt{\g^{2}+\r_{x}^{2}}}.
\end{equation}
Plugging (\ref{young}) into (\ref{M_t}), and using some properties of
$G_{\g}$ and $F_{\g}$, we get:
$$
\begin{aligned}
&M_{t}\geq \eps M_{xx}+\left( \tau
G^{'}_{\g}(\r_{x})-\frac{\r_{x}\r_{xx}}{\k_{x}^{2}}\right)M_{x}+\left(
\frac{\r_{xx}^{2}}{\k_{x}^{2}}-\frac{\r_{xxx}G^{'}_{\g}(\r_{x})}{\k_{x}}\right)M\\
&\quad\quad\quad\quad\quad\quad\quad-\frac{\g^{2}\tau^{2}}
{4\eps\sqrt{\g^{2}+\r_{x}^{2}}}-\frac{\g\g^{'}}{\sqrt{\g^{2}+\r_{x}^{2}}}.
\end{aligned}
$$

\noindent {\bf Step 2. }\textsf{(The boundary conditions for $M$)}\\

\noindent The boundary conditions (\ref{bc}), and the PDEs of system (\ref{pre_app_model}) imply the
following equalities on the boundary (using the smoothness of the solution up
to the boundary),
\begin{equation}\label{abcrk}
\left\{
\begin{aligned}
&\eps\k_{xx}+\frac{\r_{x}\r_{xx}}{\k_{x}}-\tau\r_{x}=0\quad&\mbox{on}&\quad
\partial I\times[0,T]\\
&(1+\eps)\r_{xx}-\tau\k_{x}=0\quad&\mbox{on}&\quad \partial
I\times[0,T].
\end{aligned}
\right.
\end{equation}
In particular (\ref{abcrk}) implies
\begin{equation}\label{abcM}
M_{x}=-\frac{\tau}{1+\eps}G^{'}_{\g}(\r_{x})M\quad\mbox{on}\quad\partial
I\times[0,T].
\end{equation}
To deal with the boundary condition (\ref{abcM}), we now introduce the
following change of unknown function:
$$
\o{M}(x,t)=\cosh(\beta x)M(x,t),\quad (x,t)\in \overline{I_{T}}.
$$
We calculate $\o{M}$ on the boundary of $I$ to get:
\begin{equation}\label{Mbar_bound}
\o{M}_{x}=\left( \beta\tanh(\beta
x)-\frac{\tau}{1+\eps}G^{'}_{\g}(\r_{x})\right)\o{M}\quad
\mbox{on}\quad \partial I\times [0,T].
\end{equation}
We claim that, for any fixed time $t$, it is impossible for $\o{M}$
to have a positive minimum at the boundary of $I$. Indeed we have
the following two cases:
$$\o{M} \mbox{ has a positive minimum at } x=1 \quad\Rightarrow\quad \o{M}_{x}\leq 0;$$
$$\o{M} \mbox{ has a positive minimum at } x=-1 \quad\Rightarrow\quad
\o{M}_{x}\geq 0.$$
Both cases violate the equation
(\ref{Mbar_bound}) in the case of the choice of $\beta=\beta(\eps,\tau)$
large enough, and hence the minimum of $\o{M}$ is attained inside the interval
$I$. Direct computations give:
\begin{equation}\label{pM_tbar}
\begin{aligned}
&\o{M}_{t}\geq \eps\o{M}_{xx}+\left[ \tau
G^{'}_{\g}(\r_{x})-\frac{\r_{x}\r_{xx}}{\k_{x}^{2}}-2\beta
\eps\tanh(\beta x)\right]\o{M}_{x}-
\frac{\cosh(\b x)\g^{2}\tau^{2}}{4\eps\sqrt{\g^{2}+\r_{x}^{2}}}-
\frac{\cosh(\b x)\g\g^{'}}{\sqrt{\g^{2}+\r_{x}^{2}}}\\
&+\left[\frac{\r_{xx}^{2}}{\k_{x}^{2}}-\frac{\r_{xxx}G^{'}_{\g}(\r_{x})}{\k_{x}}-\beta
\tanh(\beta x)\left(\tau
G^{'}_{\g}(\r_{x})-\frac{\r_{x}\r_{xx}}{\k_{x}^{2}}\right)+\eps\beta^{2}(2\tanh^{2}(\beta
x)-1)\right]\o{M}.
\end{aligned}
\end{equation}

\noindent {\bf Step 3. }\textsf{(The inequality satisfied by the minimum
  of $\o{M}$)}\\

\noindent Let
$$\o{m}(t)=\min_{x\in I}\o{M}(x,t).$$
Since the minimum is attained inside $I$, and since $\o{M}$ is
regular, there exists $x_0(t)\in I$ such that
$\o{m}(t)=\o{M}(x_0(t),t)$. We remark that we have:
$$\o{M}_{x}(x_0(t),t)=0,\quad \mbox{ and }\quad\o{M}_{xx}(x_0(t),t)\geq 0,$$
and hence, using (\ref{pM_tbar}), we can write down the equation
satisfied by $\o{m}$, we get (indeed in the viscosity sense):
\begin{equation}\label{min}
\begin{aligned}
&\o{m}_{t}\geq\overbrace{\left(\frac{\r_{xx}^{2}}{\k_{x}^{2}}-\frac{\r_{xxx}G^{'}_{\g}(\r_{x})}{\k_{x}}-\beta
\tanh(\beta x)\left(\tau
G^{'}_{\g}(\r_{x})-\frac{\r_{x}\r_{xx}}{\k_{x}^{2}}\right)+\eps\beta^{2}(2\tanh^{2}(\beta
x)-1)\right)}^{{\textrm{R}}}\o{m}\\
&\hspace{2cm}-\frac{\cosh(\b
  x)\g^{2}\tau^{2}}{4\eps\sqrt{\g^{2}+\r_{x}^{2}}}-\frac{\cosh(\b
  x)\g\g^{'}}{\sqrt{\g^{2}+\r_{x}^{2}}}
\quad\mbox{at}\quad x=x_{0}(t).
\end{aligned}
\end{equation}

\noindent {\bf Step 4. }\textsf{(Estimate of the term $\textrm{R}$)}\\

\noindent We turn our attention now to the term $R$ from
(\ref{min}). Using elementary identities, we get
\begin{equation}\label{I3}
R\geq
-\frac{\r_{xxx}G^{'}_{\g}(\r_{x})}{\k_{x}}-\frac{\beta^{2}\tanh^{2}(\beta
x)}{4}\frac{\r_{x}^{2}}{\k_{x}^{2}}
-\frac{\tau^{2}}{8\eps}(G^{'}_{\g}(\r_{x}))^{2}-\eps\beta^{2}.
\end{equation}
By (\ref{minko0r0x}),  we know that
$$\o{m}(0)\geq \g^{2}(0),$$
and the continuity of $\o{m}$ preserves its positivity at least for
short time. Then, as long as $\o{m}$ is positive, we have
\begin{equation}\label{2}
\k_{x}\geq\sqrt{\g^{2}+\r_{x}^{2}}.
\end{equation}
Let
$$\widetilde{c}(t)=\|\r_{xxx}(.,t)\|_{\infty,I}.$$
By using (\ref{2}) and some basic identities, inequality (\ref{I3})
implies:
\begin{equation}\label{I4}
R \geq -\frac{\widetilde{c}}{\sqrt{\g^{2}+\r_{x}^{2}}}-c_{1},\quad
c_{1}=\frac{\beta^{2}}{4}+\frac{\tau^{2}}{8\eps}+\eps\beta^{2}.
\end{equation}

\noindent {\bf Step 5. }\textsf{(The choice of $\g$ and conclusion)}\\

\noindent When $\g^{'}\leq 0$, we deduce from (\ref{min}) and (\ref{I4})
that the function $\overline{m}$ is a viscosity super-solution of:
\begin{equation}\label{cor:reg}
\o{m}_{t}= -\left(\frac{\widetilde{c}}{\sqrt{\g^{2}+\r_{x}^{2}}}+c_{1}\right) \o{m}
-\frac{c_{2}\g^{2}}{\sqrt{\g^{2}+\r_{x}^{2}}}-
\frac{\g\g^{'}}{\sqrt{\g^{2}+\r_{x}^{2}}},\quad c_{2}=\left(\frac{\t^{2}\cosh \b}{4\eps}\right).
\end{equation}
We remind the reader that $\r_{x}=\r_{x}(x_{0}(t),t)$.
Take the function $\g$ satisfying:
$$
\left\{
\begin{aligned}
& \frac{\g^{'}}{\g}\leq -(c_{0}+\widetilde{c}),\quad c_{0}=\min(c_{1},c_{2}),\\
& \g(0)=\g_{0}/2
\end{aligned}
\right.
$$
Plug $\overline{m}=\g^{2}$ into (\ref{cor:reg}), we directly deduce that $\g^{2}$ is a
viscosity sub-solution of (\ref{cor:reg}), and the result follows by
comparison. $\hfill{\Box}$
\end{section}

%%%%%%%%%%%%%%%%%%%%%%%%%%%%%%%%%%%%%%%%%%%%%%%%%%%%%%%%%%%%%%%%%
%                                                               %
%                       Section 4                               %
%                                                               %
%              Short time existence & regularity                %
%                                                               %
%%%%%%%%%%%%%%%%%%%%%%%%%%%%%%%%%%%%%%%%%%%%%%%%%%%%%%%%%%%%%%%%%

\begin{section}{Short time existence, uniqueness, and regularity}\label{sec4}
In this section, we will prove a result of short time existence,
uniqueness and regularity of a solution of problem
(\ref{pre_app_model}),  (\ref{ic}) and
(\ref{bc}).
\subsection{Short-time existence and uniqueness of a truncated system}
We denote
$$I_{a,b}:= I\times (a,a+b),\quad a,b\geq 0.$$
Fix $T_{0}\geq 0$. Consider the following system defined on
$I_{T_{0},T}$ by:
\begin{equation}\label{sh:p1:sys}
\left\{
\begin{aligned}
&\k_{t}=\eps\k_{xx}+\frac{\r_{x}\r_{xx}}{\k_{x}}-\tau\r_{x}\quad&\mbox{on}&\quad
I_{T_{0},T}
\\
&\r_{t}=(1+\eps)\r_{xx}-\tau\k_{x}\quad&\mbox{on}&\quad I_{T_{0},T},
\end{aligned}
\right.
\end{equation}
with the initial conditions:
\begin{equation}\label{sh:p1:ic}
\k(x,T_{0})=\k^{T_{0}}(x)\quad\mbox{and}\quad\r(x,T_{0})=\r^{T_{0}}(x),
\end{equation}
and the boundary conditions:
\begin{equation}\label{sh:p1:bc}
\left\{
\begin{aligned}
&\k(0,.)=0\quad\mbox{and}\quad\k(1,.)=1 \quad &\mbox{for}&\quad T_{0}<t<T_{0}+T\\
&\r(0,.)=\r(1,.)=0,\quad &\mbox{for}&\quad T_{0}<t<T_{0}+T.
\end{aligned}
\right.
\end{equation}
\begin{rem}\textit{\textbf{(The terms $p$ and $\a$).}}\label{remember}
In all what follows, and unless otherwise precised, the terms $p$
and $\alpha\in (0,1)$ are two fixed positive real numbers such that
$$p>3\quad \mbox{and}\quad \alpha=1-3/p.$$
\end{rem}
Concerning system (\ref{sh:p1:sys}), (\ref{sh:p1:ic}) and
(\ref{sh:p1:bc}), we have the following existence and
uniqueness result.
\begin{proposition}\textit{\textbf{(Short time existence and
      uniqueness).}}\label{fixed_p_T}
Let $p>3$, and $T_{0}\geq 0$. Let
$$
\r^{T_{0}},\k^{T_{0}}\in C^{\infty}(\bar{I}\times\{T_{0}\})
$$
be two given functions such that
$\r^{T_{0}}(0)=\r^{T_{0}}(1)=\k^{T_{0}}(0)=0$, and
$\k^{T_{0}}(1)=1$. Suppose furthermore that
$$
\k^{T_{0}}_{x}\geq \g_{0}\quad\mbox{on}\quad I\times\{t=T_{0}\},
$$
and
$$
\|(D^{s}_{x}\r^{T_{0}}, D^{s}_{x}\k^{T_{0}})\|_{\infty,I}\leq M_{0}\quad\mbox{on}\quad
I\times\{t=T_{0}\}, \quad s=1,2,
$$
where $\g_{0}>0$ and $M_{0}>0$ are two given positive real numbers. Then there exists
\begin{equation}\label{jin:1}
T=T^{*}=T^{*}(M_{0},\g_{0},\eps,\t,p)>0,
\end{equation}
such that the system (\ref{sh:p1:sys}), (\ref{sh:p1:ic}) and
(\ref{sh:p1:bc}) admits a unique solution
$$(\r,\k)\in (W^{2,1}_{p}(I_{T_{0},T}))^{2}.$$
Moreover, this solution satisfies
\begin{equation}\label{7n5eWl}
\k_{x}\geq \g_{0}/2\quad\mbox{on}\quad \overline{I_{T_{0},T}},
\end{equation}
and
\begin{equation}\label{sh:f5}
|\r_{x}|\leq 2M_{0}\quad\mbox{on}\quad \overline{I_{T_{0},T}}.
\end{equation}
\end{proposition}
{\bf Proof.} The short time existence is done by using a fixed point
argument. Since we are looking for solutions satisfying
(\ref{7n5eWl}) and (\ref{sh:f5}), we artificially modify
(\ref{sh:p1:sys}), and look for a solution of
\begin{equation}\label{app_app_model}
\left\{
\begin{aligned}
&\k_{t}=\eps\k_{xx}+\frac{\r_{xx}T_{2M_{0}}(\r_{x})}
{(\g_{0}/2)+(\k_{x}-\g_{0}/2)^{+}}-\tau\r_{x}\quad&\mbox{in}&\quad
I_{T_{0},T}\\
&\r_{t}=(1+\eps)\r_{xx}-\tau\k_{x}\quad&\mbox{in}&\quad
I_{T_{0},T},
\end{aligned}
\right.
\end{equation}
with the truncation function $ T_{\zeta}(x)= x
1\!\!1_{(-\zeta,\zeta)}+\zeta 1\!\!1_{\{x \geq \zeta\}} -\zeta
1\!\!1_{\{x\leq -\zeta\}}$, $\zeta>0$, and satisfying the same
initial and boundary data (\ref{sh:p1:ic}), (\ref{sh:p1:bc}). Denote
$$Y=W^{2,1}_{p}(I_{T_{0},T}).$$
For any constant $\lambda>0$, let us define $D^{\r}_{\l}$ and
$D^{\k}_{\l}$ as the two closed subsets of $Y$ given by:
$$
D^{\r}_{\l}=\{u\in Y;\;\|u_{x}\|_{p,I_{T_{0},T}}\leq \l,\;\;
u=\r^{T_{0}}
\mbox{ on }\partial^{p}I_{T_{0},T}\}
$$
and
$$
D^{\k}_{\l}=\{v\in Y;\;\|v_{x}\|_{p,I_{T_{0},T}}\leq \l,\;\;
v=\k^{T_{0}}
\mbox{ on }\partial^{p}I_{T_{0},T}\}.
$$
We choose $\lambda$ large enough such that these sets are nonempty.
Define the application $\Psi$ by:
$$
    \begin{aligned}
    \Psi:D^{\r}_{\l}\times D^{\k}_{\l}&\longmapsto& &D^{\r}_{\l}\times D^{\k}_{\l}&\\
    (\hat{\r},\hat{\k})&\longmapsto&
    &\Psi(\hat{\r},\hat{\k})=(\r,\k),&
    \end{aligned}
$$
where $(\r,\k)$ is a solution of the following system:
\begin{equation}\label{sh:sys_Theo1}
\left\{
\begin{aligned}
&\k_{t}=\eps
\k_{xx}+\frac{\r_{xx}T_{2M_{0}}(\hat{\r}_{x})}{(\g_{0}/2)+(\hat{\k}_{x}-\g_{0}/2)^{+}}-\tau
\hat{\r}_{x}\quad &\mbox{in}& \quad I_{T_{0},T},\\
& \r_{t}=(1+\eps)\r_{xx}-\tau \hat{\k}_{x}\quad &\mbox{in}& \quad
I_{T_{0},T},
\end{aligned}
\right.
\end{equation}
with the same initial and boundary conditions given by (\ref{sh:p1:ic})
and (\ref{sh:p1:bc}) respectively. The existence of the solution of
(\ref{sh:sys_Theo1}), (\ref{sh:p1:ic}) and (\ref{sh:p1:bc}) is a direct
consequence of Theorem \ref{classpara}. Taking
$\bar{\r}(x,t)=\r(x,t)-\r^{T_{0}}(x)$ and $\bar\k(x,t)=\k(x,t)-\k^{T_{0}}(x)$,
we can easily check that $(\bar{\r},\bar{\k})$ satisfies a parabolic
system similar to (\ref{sh:sys_Theo1}) with $(\bar{\r},\bar{\k})=0$ on
$\partial^{p}I_{T_{0},T}$. Using Sobolev estimates for parabolic
equations to the system satisfied by $(\bar{\r},\bar{\k})$, particularly
(\ref{parabolicinq}), we deduce that for sufficiently small $T>0$, we
have  $\|\r_{x}\|_{p,I_{T_{0},T}}\leq \l$,
$\|\k_{x}\|_{p,I_{T_{0},T}}\leq \l$, and hence the application $\Psi$ is
well defined.\\

\noindent {\bf The application $\Psi$ is a contraction map.} Let
$\Psi(\hat{\r},\hat{\k})=(\r,\k)$ and
$\Psi(\hat{\r}',\hat{\k}')=(\r',\k')$. Direct computations, using in
particular (\ref{parabolicinq}), give:
\begin{equation}\label{sh:f12}
\|\r-\r'\|_{Y}\leq c\sqrt{T}\|\hat \k-\hat{\k}'\|_{Y},
\end{equation}
and
\begin{equation}\label{sh:f8}
\|\k-\k'\|_{Y}\leq c\|F\|_{p,I_{T_{0},T}},
\end{equation}
with the function $F$ satisfying:
\begin{equation}\label{sh:est_F}
\begin{aligned}
&F+\tau(\hat{\r}-\hat{\r}')_{x}=\overbrace{\frac{T_{2M_{0}}(\hat{\r}_{x})}
{(\g_{0}/2)+(\hat{\k}_{x}-\g_{0}/2)^{+}}(\r_{xx}-\r'_{xx})}^{A_1}
+\overbrace{\frac{\r'_{xx}(T_{2M_{0}}
(\hat{\r}_{x})-T_{2M_{0}}(\hat{\r}'_{x}))}
{(\g_{0}/2)+(\hat{\k}_{x}-\g_{0}/2)^{+}}}^{A_2}
\\
& +
\overbrace{\r'_{xx}T_{2M_{0}}(\hat{\r}'_{x})\left(\frac{1}{(\g_{0}/2)+(\hat{\k}_{x}-\g_{0}/2)^{+}}
-\frac{1}{(\g_{0}/2)+(\hat{\k}'_{x}-\g_{0}/2)^{+}}\right)}^{A_3}.
\end{aligned}
\end{equation}
In order to prove the contraction for some small $T>0$, we need to
estimate all the terms appearing in (\ref{sh:est_F}). The term
$A_{1}$ can be easily handled. However, for the term $A_{2}$, we
proceed as follows. We apply the $L^{\infty}$ control of the spatial
derivative (see Lemma \ref{cor:con_hol_sob}) to the function
$\hat\r-\hat\r'$, we get:
\begin{equation}\label{sh:cor:A2}
\|(\hat\r-\hat\r')_{x}\|_{\infty,I_{T_{0},T}}\leq
cT^{\frac{p-3}{2p}}\|\hat\r-\hat\r'\|_{Y}.
\end{equation}
For the term $\r'_{xx}$, we apply (\ref{parabolicinq}), and hence we deduce that
\begin{equation}\label{sh:cor:A22}
\|\r_{xx}'\|_{p,I_{T_{0},T}}\leq c(M_{0}+\l).
\end{equation}
From (\ref{sh:cor:A2}) and (\ref{sh:cor:A22}), we deduce that
$$
\|A_{2}\|_{p,I_{T_{0},T}}\leq
c\frac{(M_{0}+\l)}{\g_{0}}T^{\frac{p-3}{2p}}\|\hat\r-\hat\r'\|_{Y}.
$$
The term $A_{3}$ could be treated in a similar way as the term
$A_{2}$. The above arguments, particularly (\ref{sh:f12}) and
(\ref{sh:f8}), give the contraction of $\Psi$ for small time \break
$T=T^{*}(M_{0},\g_{0},\eps,\t,p)>0$. Finally, inequalities
(\ref{7n5eWl}) and (\ref{sh:f5}) directly follow using the Sobolev
embedding in Hölder spaces (Lemma \ref{lemmalsu1}). $\hfill{\Box}$
\subsection{Regularity of the solution}
This subsection is devoted to show that the solution of
(\ref{sh:p1:sys}), (\ref{sh:p1:ic}) and
(\ref{sh:p1:bc}) enjoys more
regularity than the one indicated in Proposition \ref{fixed_p_T}.
This will be done using a special bootstrap argument, together with
the Hölder regularity of solutions of parabolic equations.
\begin{proposition}\textit{\textbf{(Regularity of the solution: bootstrap
      argument).}}\label{Bootstrap}
Under the same hypothesis of Proposition \ref{fixed_p_T}, let
$\r^{T_{0}}$ and $\k^{T_{0}}$ satisfy:
\begin{equation}\label{cor:comp_reg1}
\left\{
\begin{aligned}
&(1+\eps)\r^{T_{0}}_{xx}=\t\k^{T_{0}}_{x}\quad&\mbox{at}&\quad \partial I,\\
&(1+\eps)\k^{T_{0}}_{xx}=\t\r^{T_{0}}_{x}\quad&\mbox{at}&\quad \partial I.
\end{aligned}
\right.
\end{equation}
Then the unique solution $(\r,\k)$ given by Proposition
\ref{fixed_p_T} is in fact more regular. Precisely, it satisfies for
$\a=1-{3}/{p}$:
\begin{equation}\label{cor:result_regr}
\r,\k\in C^{3+\a,\frac{3+\a}{2}}(\overline{I_{T_{0},T}})\cap
C^{\infty}(\bar{I}\times (T_{0},T_{0}+T)),
\end{equation}
where $T$ is the time given by Proposition \ref{fixed_p_T}.
\end{proposition}
{\bf Proof.} For the sake of simplicity, let us suppose that
$T_{0}=0$.\\

\noindent {\bf The Hölder regularity.} Since $\k\in W^{2,1}_{p}(I_{T})$,
we use Lemma \ref{lemmalsu1}
to deduce that $\k_{x}\in C^{\a,\a/2}(\overline{I_{T}})$. We apply the
Hölder theory for parabolic equations Theorem \ref{holder_sol}, to the
second equation of (\ref{sh:p1:sys}) (using in particular the regularity
of the initial data $\r^{0}$), we deduce that:
\begin{equation}\label{sh:cor1:1}
\r\in C^{2+\a,1+\a/2}(\overline{I_{T}}).
\end{equation}
Here the compatibility condition is satisfied by
(\ref{cor:comp_reg1}). Using (\ref{sh:cor1:1}) and (\ref{7n5eWl}), we
deduce that $\frac{\r_{x}\r_{xx}}{\k_{x}}-\tau\r_{x}\in
C^{\a,\a/2}(\overline{I_{T}})$ and similar arguments as above give that:
\begin{equation}\label{sh:cor1:2}
\k\in C^{2+\a,1+\a/2}(\overline{I_{T}}).
\end{equation}
Repeating the above arguments, using this time (see (\ref{sh:cor1:2}))
that $\k_{x}\in
C^{1+\a,\frac{1+\a}{2}}(\overline{I_{T}})$, and
hence
\begin{equation}\label{sh:sh:cor1:1}
\r\in C^{3+\a,\frac{3+\a}{2}}(\overline{I_{T}}),
\end{equation}
where (\ref{sh:sh:cor1:1}) directly implies that
$\frac{\r_{x}\r_{xx}}{\k_{x}}-\tau\r_{x}\in
C^{1+\a,\frac{1+\a}{2}}(\overline{I_{T}})$, and therefore
\begin{equation}\label{sh:sh:cor1:2}
\k\in C^{3+\a,\frac{3+\a}{2}}(\overline{I_{T}}).
\end{equation}
The compatibility condition of order $1$ which is needed to apply
Theorem \ref{holder_sol} is always satisfied by
(\ref{cor:comp_reg1}). The Hölder regularity of $(\r,\k)$ directly
follows from (\ref{sh:sh:cor1:1}) and (\ref{sh:sh:cor1:2}).\\

\noindent {\bf The $C^{\infty}$ regularity.} In order to get the
$C^{\infty}$ regularity, we argue as in the case of the Hölder
regularity (bootstrap argument). In this case the compatibility
condition is replaced by multiplying by a test function that vanishes
near $t=0$. $\hfill{\Box}$.
\end{section}

%%%%%%%%%%%%%%%%%%%%%%%%%%%%%%%%%%%%%%%%%%%%%%%%%%%%%%%%%%
%                                                        %
%                  Section 5                             %
%                                                        %
%              A priori estimates                        %
%                                                        %
%%%%%%%%%%%%%%%%%%%%%%%%%%%%%%%%%%%%%%%%%%%%%%%%%%%%%%%%%%

\begin{section}{Exponential bounds}\label{sec5}
In this section, we will give some exponential bounds of the
solution given by Proposition \ref{fixed_p_T}, and having the
regularity shown by Proposition \ref{Bootstrap}. It is very
important, throughout all this section, to precise our notation
concerning the constants that may certainly vary from line to line.
Let us mention that a constant depending on time will be denoted by
$c(T)$. Those which do not depend on $T$ will be simply denoted by
$c$. In all other cases, we will follow the changing of the
constants in a precise manner.
\begin{proposition}\textit{\textbf{(Exponential bound in time for
  $\r_{x}$ and $\k_{x}$).}}\label{prop:jin:1}
Let
$$\r,\k\in
C^{3+\a,\frac{3+\a}{2}}(\bar{I}\times[0,\infty))\cap
C^{\infty}(\bar{I}\times(0,\infty)),$$
be a solution of (\ref{pre_app_model}), (\ref{ic}) and (\ref{bc}), with
$\r^{0}(0)=\r^{0}(1)=0$, $\k^{0}(0)=0$ and $\k^{0}(1)=1$.
Suppose furthermore that the function
$$
B=\frac{\r_{x}}{\k_{x}}\quad \mbox{satisfies}\quad
\|B\|_{L^{\infty}(I\times (0,\infty))}\leq 1.
$$
Then, for small $T^{*}=T^{*}(\eps,\tau,p)>0$, and
$A=1+\|\r^{0}\|_{W^{2-2/p}_{p}(I)}+\|\k^{0}\|_{W^{2-2/p}_{p}(I)}$,
we have for all $t\geq 0$:
\begin{equation}\label{VVimp1}
|\r_{x}|^{(\a)}_{I_{t,T^{*}}}, |\k_{x}|^{(\a)}_{I_{t,T^{*}}}\leq
cAe^{ct},
\end{equation}
and $c$ is a fixed constant independent of the initial data.
\end{proposition}
\noindent{\bf Proof.} We use the special coupling of the system
(\ref{pre_app_model}) to find our \textit{a priori} estimate.
Roughly speaking, the fact that $\k_{x}$ appears as a source term in
the second equation of system (\ref{pre_app_model}) permits, by the
$L^{p}$ theory for parabolic equations, to have $L^{p}$ bounds, in
terms of $\|\k_{x}\|_{p,I_{T}}$, on $\r_{x}$ and $\r_{xx}$ which in
their turn appear in the source terms of the first equation of
(\ref{pre_app_model}) satisfied by $\k$. All this permit to deduce
our estimates. To be more precise, let $T>0$ an arbitrarily fixed
time, the proof is divided into four steps:\\

\noindent {\bf Step 1. }\textsf{(estimating $\k_{x}$ in the $L^p$ norm)}\\

\noindent Let $\k^{'}$ be the solution of the following equation:
\begin{equation}\label{cc:slav}
\left\{
\begin{aligned}
&\k_{t}^{'}=\k_{xx}^{'} \quad &\mbox{on}& \quad I_{T}\\
&\k^{'}=\k \quad &\mbox{on}& \quad \partial^{p} I_{T}.
\end{aligned}
\right.
\end{equation}
As a solution of a parabolic equation, we use the $L^{p}$ parabolic
estimate (\ref{parabolic_est}) to the function $\k^{'}$ to deduce
that:
\begin{equation}\label{cc:slav1}
\|\k^{'}\|_{W^{2,1}_{p}(I_{T})}\leq
c(T)\left(\|\k^{0}\|_{W^{2-2/p}_{p}(I)}+1\right),
\end{equation}
where the term $1$ comes from the value of $\k'=\k$ on ${S}_{T}$. Take
\begin{equation}\label{cc:Afi}
\bar{\k}=\k-\k^{'},
\end{equation}
then the system satisfied by $\bar{\k}$ reads:
$$
\left\{
\begin{aligned}
&
\bar{\k}_{t}=\bar{\k}_{xx}-(\k^{'}_{t}-\eps\k^{'}_{xx})+\frac{\r_{x}\r_{xx}}{\k_{x}}-\tau\r_{x}\quad\mbox{on}\quad
I_{T}\\
& \bar{\k}=0\quad \mbox{on}\quad \partial^{p} I_{T}.
\end{aligned}
\right.
$$
Using the special version (\ref{parabolicinq}) of the parabolic $L^{p}$ estimate
 to the function $\bar{\k}$, we obtain:
\begin{equation}\label{jin:est_bark}
\|\bar{\k}_{x}\|_{p,I_{T}}\leq c\sqrt{T}\left(\|\k^{'}_{t}\|_{p,I_{T}}+\|\k^{'}_{xx}\|_{p,I_{T}}+\|\r_{xx}\|_{p,I_{T}}+
\|\r_{x}\|_{p,I_{T}}\right),
\end{equation}
where we have plugged into the constant $c$
the terms $\eps$, $\tau$, $p$ and $\|B\|_{\infty}$. Combining (\ref{cc:slav1}), (\ref{cc:Afi})
and (\ref{jin:est_bark}), we get:
\begin{equation}\label{ccc:aFi}
\|\k_{x}\|_{p,I_{T}}\leq c(T)\left(\|\k^{0}\|_{W^{2-2/p}_{p}(I)}+1
\right) +c\sqrt{T}\|\r\|_{W^{2,1}_{p}(I_{T})}.
\end{equation}
The term $\|\r\|_{W^{2,1}_{p}(I_{T})}$ appearing in the previous
inequality is going to be estimated in the
next step.\\

\noindent {\bf Step 2. }\textsf{(estimating $\r$ in the $W^{2,1}_{p}$ norm)}\\

\noindent As in Step 1, let $\r^{'}$, $\bar{\r}$ be the two
functions defined similarly as $\k^{'}$, $\bar{\k}$ respectively
(see (\ref{cc:slav}) and (\ref{cc:Afi})). The function $\r^{'}$
satisfies an inequality similar to (\ref{cc:slav1}) that reads:
\begin{equation}\label{cc:slav2}
\|\r^{'}\|_{W^{2,1}_{p}(I_{T})}\leq c(T)
\|\r^{0}\|_{W^{2-2/p}_{p}(I)}.
\end{equation}
The term $1$ disappeared here because $\r^{'}=\r=0$ on
$\overline{S_{T}}$. We write the system satisfied by $\bar{\r}$, we
obtain:
$$
\left\{
\begin{aligned}
&
\bar{\r}_{t}=(1+\eps)\bar{\r}_{xx}+
((1+\eps)\r^{'}_{xx}-\r^{'}_{t})-\tau\k_{x}\quad\mbox{on}\quad I_{T}\\
& \bar{\r}(x,0)=0\quad\mbox{on}\quad \partial^{p}I_{T},
\end{aligned}
\right.
$$
hence the following estimate on $\bar{\r}$, due to the special $L^{p}$
interior estimate (\ref{parabolicinq}), holds:
\begin{equation}\label{jin:estbarr}
\|\bar\r\|_{W^{2,1}_{p}(I_{T})}\leq
c\left(\|\r^{'}_{t}\|_{p,I_{T}}+\|\r^{'}_{xx}\|_{p,I_{T}}+\|\k_{x}\|_{p,I_{T}}\right).
\end{equation}
Again, we have plugged $\eps$, $\tau$ and $p$ into the constant $c$, and
we have assumed that $T\leq 1$. Combining (\ref{cc:slav2}) and (\ref{jin:estbarr}), we get in terms of $\r$:
\begin{equation}\label{jin:estbarr1}
\|\r\|_{W^{2,1}_{p}(I_{T})}\leq
c(T)\|\r^{0}\|_{W^{2-2/p}_{p}(I)}+ c\|\k_{x}\|_{p,I_{T}}.
\end{equation}
We will use this estimate in order to have a control on
$\|\k_{x}\|_{p,I_{T}}$ for sufficiently small time.\\

\noindent {\bf Step 3. }\textsf{(Estimate on a small time interval)}\\

\noindent From (\ref{ccc:aFi}) and (\ref{jin:estbarr1}), we
deduce that:
\begin{equation}\label{jin:mixrk}
\|\k_{x}\|_{p,I_{T}}\leq
c(T)\left(\|\k^{0}\|_{W^{2-2/p}_{p}(I)}+\|\r^{0}\|_{W^{2-2/p}_{p}(I)}+1\right)+c\sqrt{T}
\|\k_{x}\|_{p,I_{T}}.
\end{equation}
Let us remind the reader that all constants $c$ and $c(T)$ have been
changing from line to line. In fact, the important thing is whether they
depend on $T$ or not. Let
$$T^{*}=\frac{1}{2c^{2}},\quad c\mbox{ is the constant appearing in (\ref{jin:mixrk})},$$
we deduce, from (\ref{jin:mixrk}), that
$$\|\k_{x}\|_{p,I_{T^{*}}}\leq c_{3}\left(\|\k^{0}\|_{W^{2-2/p}_{p}(I)}+
\|\r^{0}\|_{W^{2-2/p}_{p}(I)}+1\right),$$
where $c_{3}=c_{3}(T^{*})>0$ is a positive constant which depends on
$T^{*}$. Recall the special coupling of
system (\ref{pre_app_model}), together with the above estimate, we can deduce that:
\begin{equation}\label{jin:expalrk}
\|(\r,\k)\|_{W^{2,1}_{p}(I_{T^{*}})}\leq  c_{4}\left(\|\k^{0}\|_{W^{2-2/p}_{p}(I)}+
\|\r^{0}\|_{W^{2-2/p}_{p}(I)}+1\right),
\end{equation}
with $c_{4}=c_{4}(T^{*})>0$ is also a positive constant depending on $T^{*}$ but
independent of the initial data.\\

\noindent {\bf Step 4. }\textsf{(The exponential estimate by iteration)}\\

\noindent Now we move to show the exponential bound. Set
$$f(t)=\|(\r,\k)\|_{W^{2,1}_{p}(I\times
  (t,t+T^{*}))},\quad\mbox{and}\quad g(t)=\|\k(\cdot,t)\|_{W^{2-2/p}_{p}(I)}+
\|\r(\cdot,t)\|_{W^{2-2/p}_{p}(I)}.$$
Using estimate (\ref{jin:expalrk}) of Lemma \ref{lemmalsu}, together
with estimate (\ref{jin:expalrk}) of Step 3, we get
$$g(T^{*}) \leq c_{5}f(0)\leq c_{5}c_{4} (g(0)+1),\quad c_{5}=c_{5}(T^{*}).$$
In this case, the Sobolev embedding in Hölder spaces (see Lemma
\ref{lemmalsu}), and the time iteration give immediately the result.
$\hfill{\Box}$

%=----------------------------------------------------------------------
%=----------------------------------------------------------------------
%=----------------------------------------------------------------------
%=----------------------------------------------------------------------

\begin{proposition}\textit{\textbf{(Exponential bound in time for
  $\r_{xx}$).}}\label{prop:jin:2}
Under the same hypothesis of Proposition \ref{prop:jin:1}, and for some
$T^{*}=T^{*}(\eps,\tau,p)>0$, we have:
\begin{equation}\label{exp_b_rxx}
|\r|^{(2+\a)}_{I_{t,T^{*}}}\leq c A e^{ct},\quad t\geq 0,
\end{equation}
where $A=1+\|\r^{0}\|_{W^{2-2/p}_{p}(I)}+\|\k^{0}\|_{W^{2-2/p}_{p}(I)}+
|\r^{0}|^{(2+\a)}_{I}$, and $c>0$ is a fixed positive constant independent of the initial data.
\end{proposition}
{\bf Proof.} The proof is very similar to the proof of Proposition
\ref{prop:jin:1}. It uses in particular the Hölder estimate for
parabolic equations (namely (\ref{cor:hol_ineq})), the Hölder
embedding in Sobolev spaces (Lemma \ref{lemmalsu1}), and finally the
iteration in time. $\hfill{\Box}$
\begin{rem}
We can not obtain, using similar arguments as in the proof of
Proposition \ref{prop:jin:2}, a similar exponential bound (\ref{exp_b_rxx}) for the
term $|\k|^{(2+\a)}_{I_{t,T^{*}}}$. This is due to the presence of the
term $1/\k_{x}$ in the equation involving $\k$.
\end{rem}

\section{An upper bound for the $W^{2,1}_{2}$ norm of  $\r_{xxx}$}\label{sec5.1}
This section is devoted to give a suitable upper bound for the
$W^{2,1}_{2}$ norm of $\r_{xxx}$. This result will be a consequence
of the control of the $W^{2,1}_{2}$ norm of $\k_{t}$ and $\k_{xx}$.
The goal is to use this upper bound in the parabolic Kozono-Taniuchi
inequality (see inequality (\ref{kte_est}) of Theorem \ref{KTRineq})
in order to control the $L^{\infty}$ norm of $\r_{xxx}$.

Let us fix $T_{1}>0$. In this section (Section 6) and in the
following section (Section 7), we will obtain some estimates on the
solution $(\r,\k)$ on the time interval $(0,T)$, with
\begin{equation}\label{hyp1}
T>T_{1}>0.
\end{equation}
In these estimates, we will precise the dependence on $T$ which
involves some constants depending on $T_{1}$ that may blow up as
$T_{1}$ goes to zero.
Consider the following hypothesis:\\

\noindent { {\textit{{{\bf {(H1)}}}}}} The function $\k_{x}$
satisfies:
$$
\k_{x}(x,t)\geq \g(t)>0,
$$
where $\g(t)$ is a positive decreasing function with $\g(0)=\g_{0}/2$,
$\g_{0}\in (0,1)$.\\

Let
$$
\mathcal{D}=I_{T},
$$
we start with the following lemma.
\begin{lem}\textit{\textbf{($W^{2,1}_{2}$ bound for $\k_{t}$ and
      $\k_{xx}$).}}\label{jin:al:lem2}
Under hypothesis (H1), and under the same hypothesis of Proposition
\ref{prop:jin:1}, we have:
$$
\|\k_{t}\|_{W^{2,1}_{2}(\mathcal{D})},
\|\k_{xx}\|_{W^{2,1}_{2}(\mathcal{D})}\leq \frac{E}{\g^{4}},
$$
where
$$\g:=\g(T)\quad \mbox{and}\quad E:=de^{d T},$$
with $d\geq 1$ is a positive constant
depending on the initial conditions but independent of $T$, and will
be given at the end of the proof.
\end{lem}
\begin{rem}\textit{\textbf{(The constant $E$ depending on time).}}\label{Go3nan1}
Let us stress on the fact that, throughout the proof and in the rest
of the paper, the term $E=de^{dT}$ of Lemma \ref{jin:al:lem2} might
vary from line to line. In other words, the term $d$ in the
expression of $E$ might certainly vary from line to line, but always
satisfying the fact of just being dependent on the initial data of
the problem. The different $E$'s appearing in different estimates
can be made the same by simply taking the maximum between them.
Therefore they will all be denoted by the same letter $E$.
\end{rem}
\noindent {\bf Proof of Lemma \ref{jin:al:lem2}.} Define the functions $u$ and $v$ by:
$$u(x,t)=\r_{t}(x,t)\quad \mbox{and}\quad v(x,t)=\k_{t}(x,t).$$
We write down the equations satisfied by $u$ and $v$ respectively:
\begin{equation}\label{7is11}
\left\{
\begin{aligned}
&u_{t}=(1+\eps)u_{xx}-\tau v_{x}\quad &\mbox{on}&\quad \mathcal{D},\\
&u|_{S_{T}}=0,\\
& u|_{t=0}=u^{0}:=(1+\eps)\r^{0}_{xx}-\tau \k^{0}_{x}&\mbox{on}&\quad I,
\end{aligned}
\right.
\end{equation}
\begin{equation}\label{7is12}
\left\{
\begin{aligned}
&v_{t}=\eps
v_{xx}+\frac{\r_{xx}}{\k_{x}}u_{x}+Bu_{xx}-B\frac{\r_{xx}}{\k_{x}}v_{x}-\tau
u_{x}\quad &\mbox{on}&\quad \mathcal{D},\\
&v|_{S_{T}}=0,\\
& v|_{t=0}=v^{0}:=\eps
\k^{0}_{xx}+\frac{\r^{0}_{x}\r^{0}_{xx}}{\k^{0}_{x}}-\tau \r^{0}_{x}&\mbox{on}&\quad I.
\end{aligned}
\right.
\end{equation}
The proof could be divided into three steps. As a first step, we will
estimate the $L^{\infty}(\mathcal{D})$ norm of the term $v_{x}=\k_{tx}$. In the
second step, we will control the $W^{2,1}_{2}(\mathcal{D})$ norm of
$v=\k_{t}$. Finally, in the third step, we will show how to deduce a
similar control on the $W^{2,1}_{2}(\mathcal{D})$ norm of $\k_{xx}$.\\

\noindent {\bf Step 1. }\textsf{(Estimating $\|\k_{xt}\|_{\infty, \mathcal{D}}$)}\\

\noindent It is worth recalling the equation
satisfied by $\k$:
$$
\k_{t}=\eps\k_{xx}+\frac{\r_{x}\r_{xx}}{\k_{x}}-\tau\r_{x}.
$$
In Proposition \ref{Bootstrap}, we have shown that $\k\in
C^{3+\a,\frac{3+\a}{2}}$. Therefore, writing the parabolic Hölder
estimate (see (\ref{cor:hol_ineq})), we obtain:
\begin{equation}\label{minui1}
\|\k_{tx}\|_{\infty, \mathcal{D}}\leq |\k|^{(3+\a)}_{\mathcal{D}}\leq
c^{H}\left(1+ \left|\frac{\r_{x}\r_{xx}}{\k_{x}}\right|^{(1+\a)}_{\mathcal{D}}+\left|\r_{x}
  \right|^{(1+\a)}_{\mathcal{D}}\right),
\end{equation}
where the term $1$ comes from the boundary conditions, and $c^{H}>0$
is a positive constant that can be estimated as $c^{H}\leq E$ (see
Remark \ref{notwritten}). We use the elementary identity
$$|fg|_{\mathcal{D}}^{(1+\a)}\leq
\|f\|_{\infty,\mathcal{D}}|g|_{\mathcal{D}}^{(1+\a)}
+\|g\|_{\infty,\mathcal{D}}|f|_{\mathcal{D}}^{(1+\a)}
+\|f_{x}\|_{\infty,\mathcal{D}}|g|_{\mathcal{D}}^{(\a)}
+\|g_{x}\|_{\infty,\mathcal{D}}|f|_{\mathcal{D}}^{(\a)},$$
to the term
$\left|\frac{\r_{x}\r_{xx}}{\k_{x}}\right|^{(1+\a)}_{\mathcal{D}}$ with
$f=\frac{\r_{x}}{\k_{x}}$ and $g=\r_{xx}$, we get:
\begin{eqnarray}\label{azar25}
\left|\frac{\r_{x}\r_{xx}}{\k_{x}}\right|^{(1+\a)}_{\mathcal{D}}&\leq&
3|\r|_{\mathcal{D}}^{(3+\a)}+\|\r_{xx}\|_{\infty,\mathcal{D}}
\left\<\frac{\r_{x}}{\k_{x}}\right\>_{\mathcal{D}}^{(1+\a)}+
\|\r_{xxx}\|_{\infty,\mathcal{D}}\left\<\frac{\r_{x}}{\k_{x}}\right\>_{\mathcal{D}}^{(\a)}\nonumber\\
&+&\frac{2|\r|_{\mathcal{D}}^{(2+\a)}}{\g}(\|\r_{xx}\|_{\infty,\mathcal{D}}+\|\k_{xx}\|_{\infty,\mathcal{D}}),
\end{eqnarray}
where we have used the fact that $\k_{x}\geq \g$ and $\k_{x}\geq
|\r_{x}|$. We plug (\ref{azar25}) in (\ref{minui1}), we obtain:
\begin{equation}\label{minui3}
\|\k_{tx}\|_{\infty, \mathcal{D}}\leq
E\left(1+|\r|_{\mathcal{D}}^{(3+\a)}
+\left\<\frac{\r_{x}}{\k_{x}}\right\>_{\mathcal{D}}^{(1+\a)}
+|\r|_{\mathcal{D}}^{(3+\a)}\left\<\frac{\r_{x}}{\k_{x}}\right\>_{\mathcal{D}}^{(\a)}
+\frac{1}{\g}\left(1+|\k|_{\mathcal{D}}^{(2+\a)}\right)\right),
\end{equation}
where we have used used the fact that the term
$|\r|^{(2+\a)}_{\mathcal{D}}$ has an exponential bound (see Proposition
\ref{prop:jin:2}) of the form $|\r|^{(2+\a)}_{\mathcal{D}}\leq E$.\\

\noindent {\bf Step 1.1. }$\left(\textsf{Estimating}
  \left\<\frac{\r_{x}}{\k_{x}}\right\>^{(1+\a)}_{\mathcal{D}}\right)$\\

\noindent From the definition of the Hölder norm, we see that in
order to control $
\left\<\frac{\r_{x}}{\k_{x}}\right\>^{(1+\a)}_{\mathcal{D}}$, it
suffices to control the three quantities:
$$\left\<\frac{\r_{x}}{\k_{x}}\right\>^{\left(\frac{1+\a}{2}\right)}_{t,\mathcal{D}},\quad
\left\<\left(\frac{\r_{x}}{\k_{x}}\right)_{x}\right\>^{(\a)}_{x,\mathcal{D}},\quad
\mbox{and} \quad
\left\<\left(\frac{\r_{x}}{\k_{x}}\right)_{x}\right\>^{\left(\frac{\a}{2}\right)}_{t,\mathcal{D}}.$$
We use the the following identity:
$$\left\<\frac{f}{g}\right\>^{\left(\a\right)}_{t,\mathcal{D}}\leq
\left\|\frac{f}{g} \right\|_{\infty,\mathcal{D}}\left\|\frac{1}{g}
\right\|_{\infty,\mathcal{D}}\left\<g\right\>^{\left(\a\right)}_{t,\mathcal{D}}+
\left\|\frac{1}{g}\right\|_{\infty,\mathcal{D}}\left\<f\right\>^{\left(\a\right)}_{t,\mathcal{D}},
$$
with $f=\r_{x}$ and $g=\k_{x}$, we get
\begin{equation}\label{BHest1}
\left\<\frac{\r_{x}}{\k_{x}}\right\>^{\left(\frac{1+\a}{2}\right)}_{t,\mathcal{D}}\leq
\frac{1}{\g}
\left(\<\r_{x}\>^{\left(\frac{1+\a}{2}\right)}_{t,\mathcal{D}}+
\<\k_{x}\>^{\left(\frac{1+\a}{2}\right)}_{t,\mathcal{D}} \right).
\end{equation}
Similarly, we obtain:
\begin{equation}\label{BHest2}
\left\<\frac{\r_{xx}}{\k_{x}}\right\>^{(\a)}_{x,\mathcal{D}}\leq
\frac{\|\r_{xx}\|_{\infty,\mathcal{D}}}{\g^{2}}
\<\k_{x}\>^{(\a)}_{x,\mathcal{D}}+\frac{\<\r_{xx}\>^{(\a)}_{x,\mathcal{D}}}{\g}.
\end{equation}
We also use the inequality:
$$\left\<fg\right\>^{(\a)}_{x,\mathcal{D}}\leq
  \|f\|_{\infty,\mathcal{D}}\<g\>_{x,\mathcal{D}}^{(\a)}
+\|g\|_{\infty,\mathcal{D}}\<f\>_{x,\mathcal{D}}^{(\a)},$$
with $f=\frac{\k_{xx}}{\k_{x}}$ and $g=\frac{\r_{x}}{\k_{x}}$, we get:
\begin{equation}\label{BHest3}
\left\<\frac{\k_{xx}\r_{x}}{\k^{2}_{x}}\right\>^{(\a)}_{x,\mathcal{D}}\leq
\frac{\<\k_{xx}\>^{(\a)}_{x,\mathcal{D}}}{\g}+
\frac{\|\k_{xx}\|_{\infty,\mathcal{D}}}{\g^{2}}\<\r_{x}\>^{(\a)}_{x,\mathcal{D}}+
\frac{\|\k_{xx}\|_{\infty,\mathcal{D}}}{\g^{2}}\<\k_{x}\>^{(\a)}_{x,\mathcal{D}}.
\end{equation}
Similarly, we get
\begin{equation}\label{BHest4}
\left\<\frac{\r_{xx}}{\k_{x}}\right\>^{\left(\frac{\a}{2}\right)}_{t,\mathcal{D}}\leq
\frac{\|\r_{xx}\|_{\infty,\mathcal{D}}}{\g^{2}}
\<\k_{x}\>^{\left(\frac{\a}{2}\right)}_{t,\mathcal{D}}+\frac{\<\r_{xx}\>^{\left(\frac{\a}{2}\right)}_{t,\mathcal{D}}}{\g},
\end{equation}
and
\begin{equation}\label{BHest5}
\left\<\frac{\k_{xx}\r_{x}}{\k^{2}_{x}}\right\>^{\left(\frac{\a}{2}\right)}_{t,\mathcal{D}}\leq
\frac{\<\k_{xx}\>^{\left(\frac{\a}{2}\right)}_{t,\mathcal{D}}}{\g}+
\frac{\|\k_{xx}\|_{\infty,\mathcal{D}}}{\g^{2}}\<\r_{x}\>^{\left(\frac{\a}{2}\right)}_{t,\mathcal{D}}+
\frac{\|\k_{xx}\|_{\infty,\mathcal{D}}}{\g^{2}}\<\k_{x}\>^{\left(\frac{\a}{2}\right)}_{t,\mathcal{D}}.
\end{equation}
Collecting the above inequalities (\ref{BHest1}), (\ref{BHest2}),
(\ref{BHest3}), (\ref{BHest4}), and (\ref{BHest5}) yield:
\begin{equation}\label{BHestsette}
\left\<\frac{\r_{x}}{\k_{x}}\right\>^{(1+\a)}_{\mathcal{D}}\leq
\frac{E}{\g^{2}}\left(1+|\k|_{\mathcal{D}}^{(2+\a)}+\|\k_{xx}\|_{\infty,\mathcal{D}}\<\k_{x}\>_{\mathcal{D}}^{(\a)}\right),
\end{equation}
where we have used the fact that $1\leq \frac{E}{\g}$, $\g\leq 1$ and
$|\r|_{\mathcal{D}}^{(2+\a)}\leq E$ (see Proposition
\ref{prop:jin:2}).\\

\noindent {\bf Step 1.2. }$\left(\textsf{Estimating }
  |\r|^{(3+\a)}_{\mathcal{D}} \textsf{ and } |\k|^{(2+\a)}_{\mathcal{D}}\right)$\\

\noindent Using Hölder estimate for parabolic equations (estimate
(\ref{cor:hol_ineq}) of Proposition \ref{holder_sol}), and similar
computations to that of the previous step, we deduce that:
\begin{equation}\label{BHest7}
|\k|^{(2+\a)}_{\mathcal{D}}\leq \frac{E}{\g} \left(1+|\k_{x}|^{(\a)}_{\mathcal{D}}\right),
\end{equation}
and
\begin{equation}\label{BHest8}
|\r|^{(3+\a)}_{\mathcal{D}}\leq \frac{E}{\g} \left(1+|\k_{x}|^{(\a)}_{\mathcal{D}}\right).
\end{equation}

\noindent {\bf Step 1.3. }\textsf{(The estimate for $\|\k_{tx}\|_{\infty,\mathcal{D}}$)}\\

\noindent By combining (\ref{minui3}), (\ref{BHestsette}),
(\ref{BHest7}), (\ref{BHest8}), and by using the fact that
$|\k_{x}|_{\mathcal{D}}^{(\a)}$ has an exponential estimate (see
estimate (\ref{VVimp1}) of Proposition \ref{prop:jin:1}), we deduce that:
\begin{equation}\label{AddRe}
|\k|^{(3+\a)}_{\mathcal{D}}\leq \frac{E}{\g^{3}},
\end{equation}
which will be useful later, and as a particular subcase, we have:
\begin{equation}\label{BHest9}
\|\k_{tx}\|_{\infty, \mathcal{D}}\leq \frac{E}{\g^{3}},
\end{equation}
where we have frequently used that $\g\leq 1$, and we have always
taken the maximum of all the exponential bounds of the $E=de^{d T}$
form.\\

\noindent {\bf Step 2. }\textsf{(Estimating $\|\k_{t}\|_{W^{2,1}_{2}(\mathcal{D})}$)}\\

\noindent {\bf Step 2.1. }\textsf{(Estimating $\|u\|_{W^{2,1}_{2}(\mathcal{D})}$)}\\

\noindent We use the $L^{2}$ estimates for parabolic equations
(Theorem \ref{classpara}) to the function $u$ satisfying
(\ref{7is11}), we obtain:
\begin{equation}\label{BHest10}
\|u\|_{W^{2,1}_{2}(\mathcal{D})}\leq E(1+\|v_{x}\|_{2,\mathcal{D}}).
\end{equation}
The term $1$ in (\ref{BHest10}) comes from estimating the initial
data $u^{0}$. Since $v_{x}=\k_{tx}$, we plug the estimate
(\ref{BHest9}) obtained in Step 1.3 into (\ref{BHest10}), we get
\begin{equation}\label{BHest11}
\|u\|_{W^{2,1}_{2}(\mathcal{D})}\leq \frac{E}{\g^{3}}.
\end{equation}

\noindent {\bf Step 2.2. }\textsf{(Estimating $\|v\|_{W^{2,1}_{2}(\mathcal{D})}$)}\\

\noindent Arguing in a similar manner as in the previous step, we
obtain the following estimate for the function $v$, the solution of
the parabolic equation (\ref{7is12}):
\begin{eqnarray}\label{BHest12}
&\hspace{-3cm}\|v\|_{W^{2,1}_{2}(\mathcal{D})}\leq
E\left(1+\left\|\frac{\r_{xx}}{\k_{x}}\right\|_{\infty,\mathcal{D}}\|u_{x}\|_{2,\mathcal{D}}
+\|B\|_{\infty,\mathcal{D}}\|u_{xx}\|_{2,\mathcal{D}}\right.\nonumber\\
&\hspace{1cm}\left.+\|B\|_{\infty,\mathcal{D}}
\left\|\frac{\r_{xx}}{\k_{x}}\right\|_{\infty,\mathcal{D}}\|v_{x}\|_{2,\mathcal{D}}+\|u_{x}\|_{2,\mathcal{D}}\right),
\end{eqnarray}
and hence, from (\ref{BHest9}), (\ref{BHest11}), and doing some
computations, we deduce from (\ref{BHest12}) that:
\begin{equation}\label{BHest13}
\|v\|_{W^{2,1}_{2}(\mathcal{D})}\leq \frac{E}{\g^{4}}.
\end{equation}
The goal of Step 2 follows since $v=\k_{t}$.\\

\noindent {\bf Step 3. }\textsf{(Estimating
$\|\k_{xx}\|_{W^{2,1}_{2}(\mathcal{D})}$)}\\

\noindent The estimate of $\|\k_{xx}\|_{W^{2,1}_{2}(\mathcal{D})}$
requires a special attention. We will mainly use the equations on
$\r$ and $\k$. The four parts $\|\k_{xx}\|_{2,\mathcal{D}}$,
$\|\k_{xxt}\|_{2,\mathcal{D}}$, $\|\k_{xxx}\|_{2,\mathcal{D}}$ and
$\|\k_{xxxx}\|_{2,\mathcal{D}}$ of
the above norm will be estimated separately.\\

\noindent {\bf Step 3.1. }\textsf{(Estimate of
  $\|\k_{xx}\|_{2,\mathcal{D}}$)}\\

\noindent Inequality (\ref{BHest7}) directly implies that

\begin{equation}\label{Ayst2150708}
\|\k_{xx}\|_{\infty,\mathcal{D}}\leq \frac{E}{\g},
\end{equation}
hence $\|\k_{xx}\|_{2,\mathcal{D}}\leq \frac{E}{\g}$.\\

\noindent {\bf Step 3.2. }\textsf{(Estimate of
  $\|\k_{xxxx}\|_{2,\mathcal{D}}$)}\\

\noindent We first derive the equation on $\r$ two times in $x$, we
deduce (using (\ref{BHest11})) that $\|\r_{xxxx}\|_{2,\mathcal{D}}$
has the same upper bound as $\|\k_{xxx}\|_{2,\mathcal{D}}$, i.e.
\begin{equation}\label{korz:10}
\|\r_{xxxx}\|_{2,\mathcal{D}}\leq \frac{E}{\g^{3}}.
\end{equation}
We derive the equation on $\k$ two times with respect to the
variable $x$, we obtain:
\begin{eqnarray*}
\k_{txx}&=&\eps \k_{xxxx}+\frac{2\r_{xx}\r_{xxx}}{\k_{x}}
-\frac{\k_{xx}\r^{2}_{xx}}{\k^{2}_{x}}+\frac{\r_{x}\r_{xxxx}}{\k_{x}}-\frac{\r_{x}\r_{xxx}\k_{xx}}{\k^{2}_{x}}\\
&-&\frac{\r^{2}_{xx}\k_{xx}}{\k^{2}_{x}}
-\frac{\r_{x}\r_{xx}\k_{xxx}}{\k^{2}_{x}}-\frac{\r_{x}\k_{xx}\r_{xxx}}{\k^{2}_{x}}+
\frac{2\k^{2}_{xx}\r_{x}\r_{xx}}{\k^{3}_{x}}-\tau \r_{xxx},
\end{eqnarray*}
and we use (\ref{korz:10}) and our controls obtained in the previous
steps, in order to deduce that:
$$
\|\k_{xxxx}\|_{2,\mathcal{D}}\leq \frac{E}{\g^{4}}.
$$
In fact, the highest power comes from estimating the following term:
$$\left\|\frac{\k_{xx}^{2}\r_{x}\r_{xx}}{\k^{3}_{x}}\right\|_{2,\mathcal{D}}\leq
\left\|\frac{\k_{xx}^{2}\r_{xx}}{\k^{2}_{x}}\right\|_{\infty,\mathcal{D}}\sqrt{T}\leq
\frac{E}{\g^{4}},$$ where we have used the $L^{\infty}$ estimate of
$\|\k_{xx}\|_{\infty,\mathcal{D}}$. All other estimates are easily
deduced. Let us just state how to estimate the other term were
$\|\k_{xx}\|_{\infty,\mathcal{D}}$ interferes. In fact, we have:
$$\left\|\frac{\r_{x}\r_{xxx}\k_{xx}}{\k^{2}_{x}}\right\|_{2,\mathcal{D}}\leq
\left\|\frac{\k_{xx}}{\k_{x}}\right\|_{\infty,\mathcal{D}}\|\r_{xxx}\|_{2,\mathcal{D}}\leq
\frac{E}{\g^{3}}.$$

\noindent {\bf Step 3.3. }\textsf{(Estimate of
  $\|\k_{xxt}\|_{2,\mathcal{D}}$ and $\|\k_{xxx}\|_{2,\mathcal{D}}$)}\\

\noindent As an immediate consequence of (\ref{BHest13}), we get
$$\|\k_{xxt}\|_{2,\mathcal{D}}\leq \frac{E}{\g^{4}}.$$
Deriving the equation on $\k$ with respect to $x$, we obtain:
\begin{equation}\label{korz:8.5}
\k_{tx}=\eps
\k_{xxx}+\frac{\r^{2}_{xx}}{\k_{x}}+\frac{\r_{x}\r_{xxx}}{\k_{x}}-
\frac{\r_{x}\k_{xx}\r_{xx}}{\k^{2}_{x}}-\tau \r_{xx}.
\end{equation}
The estimate (\ref{BHest8}) gives
$$
\|\r_{xxx}\|_{2,\mathcal{D}}\leq \frac{E}{\g},
$$
which, together with (\ref{Ayst2150708}) and (\ref{korz:8.5}), give
$\|\k_{xxx}\|_{2,\mathcal{D}}\leq \frac{E}{\g^{3}}$. We deduce as a
conclusion that:
$$
 \|\k_{xx}\|_{W^{2,1}_{2}(\mathcal{D})}\leq
\frac{E}{\g^{4}},
$$
and this terminates the proof. $\hfill{\Box}$\\

We move now to the main result of this section.
\begin{lem}\textit{\textbf{($W^{2,1}_{2}$ bound for $\r_{xxx}$ )}}\label{min:lem2}\\
Under the same hypothesis of Lemma \ref{jin:al:lem2}, we have:
\begin{equation}\label{min:eqn2}
\|\r_{xxx}\|_{W^{2,1}_{2}(\mathcal{D})}\leq \frac{E}{\g^{4}}.
\end{equation}
\end{lem}
{\bf Proof.} Set
$$\bar{\k}=\frac{\tau}{(1+\eps)^{2}}\k_{t}+\frac{\tau}{1+\eps}\k_{xx}$$ and
$$w=\r_{xxx}-\bar{\k}.$$
We write down, after doing some computations, the equation satisfied
by $w$:
\begin{equation}\label{amal:2}
\left\{
\begin{aligned}
& w_{t}=(1+\eps)w_{xx}- \frac{\tau}{(1+\eps)^{2}}\k_{tt}\quad
\mbox{on} \quad
\mathcal{D}\\
& w_{x}|_{S_{T}}=0 \quad\mbox{on}\quad S_{T}\\
& w|_{t=0}:=w^{0}=\r^{0}_{xxx}-
\frac{\tau(1+2\eps)}{(1+\eps)^{2}}\k^{0}_{xx}-
\frac{\tau}{(1+\eps)^{2}}\frac{\r^{0}_{x}\r^{0}_{xx}}{\k^{0}_{x}}+
\frac{\tau^{2}}{(1+\eps)^{2}}\r^{0}_{x}.
\end{aligned}
\right.
\end{equation}
Here $w_{x}|_{S_{T}}=0$ can be checked by deriving the equation
satisfied by $\r$ with respect to $x$ and then with respect to $t$,
and by using the equality $(1+\eps)\r_{xx}=\tau \k_{x}$ satisfied on
the boundary $\partial I$ (which is a consequence of the
compatibility conditions). Applying the $L^{2}$ theory with Neumann
conditions (see for instance \cite[Chapter 4, Section 10]{LSU}) to
(\ref{amal:2}), we get that $\r_{xxx}=w+\bar{\k}$ satisfies
\begin{equation}\label{amal:5}
\|\r_{xxx}\|_{W^{2,1}_{2}(\mathcal{D})}\leq
E\left(1+\|\k_{tt}\|_{2,\mathcal{D}}+
\|\k_{t}\|_{W^{2,1}_{2}(\mathcal{D})}+\|\k_{xx}\|_{W^{2,1}_{2}(\mathcal{D})}\right),
\end{equation}
and eventually (\ref{amal:5}) with Lemma \ref{jin:al:lem2} gives
immediately the result. $\hfill{\Box}$

%                          H       E       R       E

\section{An upper bound for the $BMO$ norm of $\r_{xxx}$}\label{sec5.2}
 This section is devoted to give a suitable upper bound for the
$BMO$ norm of $\r_{xxx}$. This result will be a consequence of the
control of the $BMO$ norm of a suitable extension of $\k_{xx}$. The
goal is to use this upper bound in the Kozono-Taniuchi inequality
(see inequality (\ref{kte_est}) of Theorem \ref{KTRineq}) in order
to control the $L^{\infty}$ norm of $\r_{xxx}$. We first give some
useful definitions.
\begin{definition}\textit{\textbf{(The ``symmetric and periodic'' extension of a
    function).}} Let $f\in C(\overline{I_{T}})$, we define $f^{sym}$ (constructed out
of $f$) over $\R\times (0,T)$, first by the symmetry of $f$ with
respect to the line $x=0$ over the interval $(-1,0)$, and then by
spatial periodicity.
\end{definition}

\begin{definition}\label{ch3:defasym}\textit{\textbf{(The ``antisymmetric and
      periodic'' extension of a
    function).}} We define the function $f^{asym}$ in a similar
    manner as $f^{sym}$, where we take the antisymmetry of $f$
    instead of the symmetry.
\end{definition}
We start with the following lemma that reflects a useful relation between the
$BMO$ norm of $f^{sym}$ and $f^{asym}$.
\begin{lem}\textit{\textbf{(A relation between $f^{sym}$ and
      $f^{asym}$).}}\label{3n2bnkzBTA} Let $f\in C(\overline{I_{T}})$, then:
$$
\|f^{sym}\|_{BMO(\R\times (0,T))}\leq c \left(\|f^{asym}\|_{BMO(\R\times
    (0,T))}+m_{2I\times (0,T)}\left(|f^{sym}|\right)\right),
$$
where $c>0$ is a universal constant.
\end{lem}
The proof of this lemma will be presented in Appendix B.
The next lemma gives a control of the $BMO$ norm of $(\k_{xx})^{asym}$.
\begin{lem}\textit{\textbf{($BMO$ bound for $(\k_{xx})^{asym}$).}}\label{jin:al:lem1}
Under hypothesis (H1), and under the same hypothesis of Proposition
\ref{prop:jin:1}, we have:
\begin{equation}\label{jin:al:1}
\|(\k_{xx})^{asym}\|_{BMO(\R\times (0,T))}\leq c e^{c T},
\end{equation}
where $c>0$ is a constant depending on the initial conditions (but
independent of $T$). The function $(\k_{xx})^{asym}$ is given via
Definition \ref{ch3:defasym}.
\end{lem}
{\bf Proof.} Let $\bar{\k}(x,t)=\k(x,t)-\k^{0}(x)$. We notice that
$\bar{\k}|_{S_{T}}=0$, therefore $\bar{\k}^{asym}$ satisfies:
\begin{equation}\label{anti:1}
\left\{
\begin{aligned}
&\bar{\k}^{asym}_{t}-\eps \bar{\k}^{asym}_{xx}=
\frac{(\r_{x})^{asym}\r^{asym}_{xx}}{(\k_{x})^{asym}}-\tau
(\r_{x})^{asym}+\eps (\k^{0}_{xx})^{asym}\quad \mbox{on} \quad
\R\times
(0,T)\\
& \bar{\k}^{asym}(x,0)=0.
\end{aligned}
\right.
\end{equation}
We already know that the right hand side of (\ref{anti:1}) is
bounded in $L^{\infty}$ by $E=ce^{cT}$, and hence (using Theorem
\ref{bmo_theory}) the result follows. $\hfill{\Box}$\\

We now present the principal result of this section.
\begin{lem}\textit{\textbf{($BMO$ bound for $\r_{xxx}$).}}\label{min:lem1}
Under the same hypothesis of Lemma \ref{jin:al:lem1}, we have:
\begin{equation}\label{min:eqn1}
\|\r_{xxx}\|_{BMO(\mathcal{D})}\leq E.
\end{equation}
\end{lem}
{\bf Proof.} Take
$$v=\r_{x}-\frac{\tau \k}{1+\eps},$$
the equation satisfied by $v$ reads:
$$
\left\{
\begin{aligned}
& v_{t}=(1+\eps)v_{xx}-\frac{\eps\tau}{1+\eps}\k_{xx}-
\frac{\tau}{1+\eps}\frac{\r_{x}\r_{xx}}{\k_{x}}+\frac{\tau^{2}}{1+\eps}
\r_{x} \quad \mbox{on} \quad \mathcal{D}\\
& v|_{t=0}=v^{0}:=\r^{0}_{x}-\frac{\tau}{1+\eps}\k^{0}\quad \mbox{on} \quad I\\
& v_{x}|_{S_{T}}=0,
\end{aligned}
\right.
$$
where we have used the compatibility conditions to check that
$v_{x}|_{S_{T}}=0$. We can assume, without loss of generality, that
the initial condition $v^{0}=0$. This is because being non-zero just
adds a constant depending on the initial conditions in the final
estimate that we are looking for. From the fact that
$v_{x}|_{S_{T}}=0$, we can easily deduce that the function $v^{sym}$
satisfies:
$$
\left\{
\begin{aligned}
&
v^{sym}_{t}=(1+\eps)v^{sym}_{xx}+\overbrace{\frac{\tau^{2}}{1+\eps}
(\r_{x})^{sym}-
\frac{\tau}{1+\eps}\frac{(\r_{x})^{sym}(\r_{xx})^{sym}}{(\k_{x})^{sym}}-\frac{\eps\tau}{1+\eps}(\k_{xx})^{sym}}^{g}\\
&\quad\quad\quad\quad\quad\quad\quad\quad\quad\quad\quad\quad\quad\quad\quad
\hspace{4.5cm}\mbox{on} \quad \R\times (0,T)\\
& v^{sym}(x,0)=0 \quad \mbox{on} \quad \R,
\end{aligned}
\right.
$$
therefore, using the $BMO$ estimate (\ref{cau:2}) for parabolic
equations, to the function $v$, one gets:
$$
\|v^{sym}_{xx}\|_{BMO(\R\times(0,T))}\leq c\left[\|g\|_{BMO(\R\times
(0,T))}+m_{2I\times (0,T)}(|g|)\right],
$$
where $c=c(T_{1})>0$ with $0<T_{1}\leq T$. From Propositions
\ref{prop:jin:1}, \ref{prop:jin:2}, we deduce that
$$
\|g\|_{BMO(\R\times (0,T))}\leq E+\|(\k_{xx})^{sym}\|_{BMO(\R\times
(0,T))},
$$
and
$$
m_{2I\times (0,T)}(|g|)\leq E+m_{2I\times (0,T)}(|(\k_{xx})^{sym}|).
$$
Recall the definition of the term $E$ from Remark \ref{Go3nan1}.
At this stage, we write the following estimate:
\begin{equation}\label{tashma}
\|(\k_{xx})^{sym}\|_{BMO(\R\times (0,T))}\leq
c\left[\|(\k_{xx})^{asym}\|_{BMO(\R\times (0,T))}+m_{2I\times
    (0,T)}(|(\k_{xx})^{sym}|)\right],
\end{equation}
which can be deduced using Lemma \ref{3n2bnkzBTA}. The constant
$c>0$ appearing in (\ref{tashma}) is independent of $T$. Finally, we
deduce that:
$$
\|v^{sym}_{xx}\|_{BMO(\R\times (0,T))} \leq
c\left(E+T^{-1/p}\|\k_{xx}\|_{p,\mathcal{D}}\right).
$$
From (\ref{hyp1}), (\ref{jin:expalrk}), we know that
$T^{-1/p}\|\k_{xx}\|_{p,\mathcal{D}}\leq T_{1}^{-1/p} E$. From the
previous two inequalities, and since $v_{xx}=\r_{xxx}-\frac{\tau
\k_{xx}}{1+\eps}$, we easily arrive to our result.  $\hfill{\Box}$

%==---------------------------------------------------------------------------------------------------
%  HERE IN THE CORRECTION
%==---------------------------------------------------------------------------------------------------

\section{$L^{\infty}$ bound for $\r_{xxx}$}\label{sec5.3}
In this section, we use the results of Sections \ref{sec5},
\ref{sec5.1} and \ref{sec5.2}, in order to give an $L^{\infty}$
bound for $\r_{xxx}$ via the Kozono-Taniuchi inequality.
\begin{proposition}\textit{\textbf{($L^{\infty}$ bound for
    $\r_{xxx}$).}}\label{bmosobrxxx}
Under hypothesis (H1), and under the same hypothesis of Proposition
\ref{prop:jin:1}, we have $\forall T>0$:
\begin{equation}\label{bmosobk_xx}
\|\r_{xxx}\|_{\infty,I_{T}}\leq
E\left(1+\log^{+}\frac{1}{\g}\right).
\end{equation}
\end{proposition}
\noindent {\bf Proof.} For $T<T_{1}=T^{*}$, where $T^{*}$ is the
short time existence result (see Theorem \ref{fixed_p_T}) given by
(\ref{jin:1}), inequality (\ref{bmosobk_xx}) directly follows. In
the other case where $T\geq T_{1}$, we apply the parabolic
Kozono-Taniuchi estimate (\ref{kte_est}) to the function $\r_{xxx}$,
together with (\ref{min:eqn1}) and (\ref{min:eqn2}). Remark that the
$\|\r_{xxx}\|_{L^{1}(I_{T})}$ can be easily estimated by the term
$E$. $\hfill{\Box}$ 
\begin{proposition}\textit{\textbf{(A priori estimates).}}\label{ven:1} Under the same hypothesis of
Proposition \ref{max_p}, the solution $(\r,\k)\in
C^{3+\a,\frac{3+\a}{2}}(\overline{I_{T}})$ satisfies for every
$0\leq t\leq T$:
\begin{equation}\label{ynumber1}
\k_{x}(.,t)\geq e^{-e^{e^{b(t+1)}}}>0,
\end{equation}
\begin{equation}\label{ven:est1}
|\r(.,t)|^{(3+\a)}_{I}\leq e^{ e^{ e^{b (t+1)}}}\quad
\mbox{and}\quad |\k(.,t)|^{(3+\a)}_{I}\leq  e^{ e^{ e^{b (t+1)}}}.
\end{equation}
Here $b>0$ is a positive constant depending on the initial conditions
and the fixed terms of the problem, but independent of time.
\end{proposition}
{\bf Proof.} Remark that if we consider a function $\g$ satisfying
(\ref{from_cras}), then the right hand side of (\ref{from_cras}) can
be estimated using (\ref{bmosobk_xx}) as follows:
\begin{equation}\label{hra12}
-(c_{0}+\|\r_{xxx}(.,t)\|_{L^{\infty}(I)})\geq -E(1+|\log
\g(t)|),\quad E=E(T)=de^{dT}.
\end{equation}
This is the motivation to consider the solution $\g_{T}$ of the
following ordinary differential equation:
\begin{equation}\label{la32}
\left\{
\begin{aligned}
& \g^{'}_{T}=-E (1 + |\log \g_{T}|)\g_{T},\quad t\in (0,T)\\
& \g_{T}(0)=\g_{0}/2,
\end{aligned}
\right.
\end{equation}
where $\g_{0}$ is given by (\ref{minko0r0x}). Then, using a
continuity argument, joint to the fact that $\overline{m}(t)\geq
\g^{2}(t)$ (see Proposition \ref{max_p}), it is easy to check that
both (\ref{hra12}) and (\ref{from_cras}) are satisfied with
$\g=\g_{T}$. Let us now define
$$\tilde{\g}(T):=\g_{T}(T)\leq \g_{T}(t)\quad \mbox{for}\quad t\in [0,T].$$
Then we deduce from (\ref{hra123}) that
$$\k_{x}(.,T)\geq \sqrt{\tilde{\g}^{2}(T)+(\r_{x}(.,T))^{2}}\geq \tilde{\g}(T),\quad \forall
T>0,$$ and finally, solving (\ref{la32}) explicitly, inequality
(\ref{ynumber1}) directly follows. Therefore, from (\ref{BHest8})
and (\ref{AddRe}), we easily deduce (\ref{ven:est1}). $\hfill{\Box}$
\end{section}

%%%%%%%%%%%%%%%%%%%%%%%%%%%%%%%%%%%%%%%%%%%%%%%%%%%%%%%%%%%
%                                                         %
%                   Section 6                             %
%                                                         %
%             Long time existence                         %
%                                                         %
%%%%%%%%%%%%%%%%%%%%%%%%%%%%%%%%%%%%%%%%%%%%%%%%%%%%%%%%%%%

\begin{section}{Long time existence and uniqueness}\label{sec6}

Now we are ready to show the main result of this paper, namely
Theorem \ref{mainresult0}.\\

\noindent {\bf Proof of Theorem \ref{mainresult0}.} Define the set $\mathcal{B}$ by:
$$
\mathcal{B}=\left\{
\begin{aligned}
&T>0; \;\exists \,! \;\mbox{ solution } (\r,\k)
  \in C^{3+\a,\frac{3+\a}{2}}(\overline{I_{T}}) \mbox{ of }\\
& (\ref{pre_app_model}), (\ref{ic}) \mbox{ and } (\ref{bc}), \mbox{
    satisfying } (\ref{MT5})
\end{aligned}
\right\}.
$$
This set is non empty by the short time existence result (Theorem
\ref{fixed_p_T}). Set
$$T_{\infty}=\sup \mathcal{B}.$$
We claim that $T_{\infty}=\infty$. Assume, by contradiction that
$T_{\infty}<\infty$. In this case, let $\delta>0$ be an arbitrary
small positive constant, and apply the short time existence result
(Theorem \ref{fixed_p_T}) with $T_{0}=T_{\infty}-\delta$. Indeed, by
the tri-exponential bounds (\ref{ynumber1}) and (\ref{ven:est1}), we
deduce that the time of existence $T^{*}$ given by (\ref{jin:1}) is
in fact independent of $\delta$. Hence, choosing $\delta$ small
enough, we obtain $T_{0}+T^{*}\in \mathcal{B}$ with
$T_{0}+T^{*}>T_{\infty}$ and hence a contradiction. $\hfill{\Box}$

\end{section}

%%%%%%%%%%%%%%%%%%%%%%%%%%%%%%%%%%%%%%%%%%%%%%%%
%                                              %
%               Appendix A                     %
%                                              %
%%%%%%%%%%%%%%%%%%%%%%%%%%%%%%%%%%%%%%%%%%%%%%%%

\begin{section}{Appendix A: miscellaneous parabolic estimates}\label{secA}

\noindent {\bf A1. Proof of Lemma \ref{prec_estpara} ($L^{p}$ estimate for
  parabolic equations)}\\

\noindent As a first step, we will prove the result in the case where $\eps=1$,
and in a second step, we will move to the case $\eps>0$. It is worth
noticing that the term $c$ may take several values only depending on $p$.\\

\noindent {\bf Step 1. }\textsf{(The estimate: case $\eps=1$)}\\

\noindent Suppose $\eps=1$. Since $u=0$ on $\partial I\times [0,T]$,
we take $\tilde{u}=u^{asym}$ (see Definition \ref{ch3:defasym}).
Also consider the function $\tilde{f}=f^{asym}$. Define ${\bar{u}}$
by
$$
{\bar{u}}=\tilde{u}\phi^{n},
$$
with
$$
\left\{
\begin{aligned}
&\phi^{n}(x)=1&\mbox{if }& \;x\in (0,2n) \\
&\phi^{n}(x)=0&\mbox{if }& \;x\geq 2n+1\mbox{ or } x\leq -1.
\end{aligned}
\right.
$$
This function satisfies
$$
\left\{
\begin{aligned}
& {\bar{u}}_{t}={\bar{u}}_{xx}+{\bar{f}},
\quad&\mbox{on}&\quad \R\times (0,T)\\
& {\bar{u}}(x,0)=0,\quad&\mbox{on}&\quad \R,
\end{aligned}
\right.
$$
with
$$
\bar{f}=\tilde{f}\phi^{n}-\tilde{u}\phi^{n}_{xx}-2\tilde{u}_{x}\phi^{n}_{x}.
$$
The proof that
\begin{equation}\label{syra}
\|u_{t}\|_{p,I_{T}}+\|u_{xx}\|_{p,I_{T}}\leq c \|f\|_{p,I_{T}}
\end{equation}
can be easily deduced by applying the Calderon-Zygmund estimates to
the function $\bar{u}$ satisfying the above equation, and passing to
the limit $n \rightarrow \infty$. Now, since $u\in
W^{2,1}_{p}(I_{T})$ with $u|_{t=0}=0$, we use \cite[Lemma 4.5, page
305]{LSU} to get
\begin{equation}\label{syra1}
\|u\|_{p,I_{T}}\leq cT(\|{u}_{t}\|_{p,I_{T}}+\|{u}_{xx}\|_{p,I_{T}})
\end{equation}
and
\begin{equation}\label{syra2}
\|u_x\|_{p,I_{T}}\leq
c\sqrt{T}(\|{u}_{t}\|_{p,I_{T}}+\|{u}_{xx}\|_{p,I_{T}}).
\end{equation}
Combining (\ref{syra}), (\ref{syra1}) and (\ref{syra2}), we deduce
that
$$
\frac{1}{T}\|u\|_{p,I_{T}}+\frac{1}{\sqrt{T}}\|u_{x}\|_{p,I_{T}}+\|u_{xx}\|_{p,I_{T}}+\|u_{t}\|_{p,I_{T}}\leq
c \|f\|_{p,I_{T}}.
$$

\noindent {\bf Step 2. }\textsf{(The estimate: general case $\eps>0$)}\\

\noindent To get the general inequality, we consider the following
rescaling of the function $u$:
$$
\hat{u}(x,t)=u(x,t/\eps),\quad (x,t)\in I_{\eps T},
$$
which allows to get the desired result. $\hfill{\Box}$\\

\noindent {\bf A2. Proof of Lemma \ref{cor:con_hol_sob} ($L^{\infty}$
  control of the spatial derivative)}\\

\noindent Since $u\in W^{2,1}_{p}(I_{T})$ for $p>3$, we know from Lemma \ref{lemmalsu1} that $u_{x}\in
C^{\alpha,\alpha/2}(\overline{I_{T}})$ for $\a=1-\frac{3}{p}$. In this case,
we use the estimate (\ref{cor:estH1}) with $\delta=\sqrt{T}$, we obtain
\begin{equation}\label{cor:eln}
\|u_{x}\|_{\infty,I_{T}}\leq
c(p)\{{T}^{\frac{\a}{2}}(\|u_{t}\|_{p,I_{T}}+\|u_{xx}\|_{p,I_{T}})
+{T}^{\frac{\a}{2}-1}\|u\|_{p,I_{T}}\}.
\end{equation}
Remark that the fact that $u=0$ on the parabolic boundary $\partial
^{p}I_{T}$, and that it obviously satisfies the equation:
$$
\left\{
\begin{aligned}
&u_{t}=u_{xx}+f,\quad \mbox{with}\quad f=u_{t}-u_{xx}\\
&u=0\quad \mbox{on} \quad \partial^{p}I_{T},
\end{aligned}
\right.
$$
then we can apply estimate (\ref{parabolicinq}) to bound the term
$\|u\|_{p,I_{T}}$. Hence (\ref{cor:eln}) becomes (with a different
constant $c(p)$):
\begin{eqnarray*}
\|u_{x}\|_{\infty,I_{T}}&\leq&c(p)\{{T}^{\frac{\a}{2}}\|u_{t}-u_{xx}\|_{p,I_{T}}
+{T}^{\frac{\a}{2}-1}T\|u_{t}-u_{xx}\|_{p,I_{T}}\}\\
&\leq& c(p)T^{\frac{\a}{2}}\|u\|_{W^{2,1}_{p}(I_{T})}\\
&\leq& c(p)T^{\frac{p-3}{2p}}\|u\|_{W^{2,1}_{p}(I_{T})},
\end{eqnarray*}
and the result follows. $\hfill{\Box}$\\

\end{section}

%%%%%%%%%%%%%%%%%%%%%%%%%%%%%%%%%%%%%%%%%%%%%%%%%%%%%%%%%
%                                                       %
%                     Appendix B                        %
%                                                       %
%%%%%%%%%%%%%%%%%%%%%%%%%%%%%%%%%%%%%%%%%%%%%%%%%%%%%%%%%

\begin{section}{Appendix B: parabolic $BMO$ theory}\label{secB}
\noindent {\bf B1. Proof of Theorem \ref{bmo_theory} (A $BMO$
  estimate in the periodic case)}\\

\noindent Let $f$ be a bounded function defined on $\R\times (0,T)$
satisfying $f(x+2,t)=f(x,t)$. We extend the function $f$ to
$\R\times \R_{+}$, first by symmetry with respect to the line
$\{t=T\}$ and after that by time periodicity of period $2T$. Call
this function $\tilde{f}$. Set $\bar{u}$ as the solution of the
following equation:
\begin{equation}\label{maison1}
\left\{
\begin{aligned}
& \bar{u}_{t}=\eps\bar{u}_{xx}+\tilde{f}\quad \mbox{on} \quad
\R\times \R_{+}\\
& \bar{u}(x,0)=0.
\end{aligned}
\right.
\end{equation}
We apply the standard result of  $BMO$ theory for parabolic
equations. Since \break $f\in L^{\infty}(\R\times (0,T))$, then
$\tilde{f}\in BMO(\R\times\R_{+})$, and hence we obtain \break that
$\bar{u}_{t}, \bar{u}_{xx}\in BMO(\R\times\R_{+})$, with the
following estimate:
\begin{equation}\label{ch3:maison-1}
\|\bar{u}_{t}\|_{BMO(\R\times\R_{+})}+\|\bar{u}_{xx}\|_{BMO(\R\times\R_{+})}\leq
c \|\tilde{f}\|_{BMO(\R\times\R_{+})},
\end{equation}
and hence (from the definition of the $BMO$ space),
\begin{equation}\label{maison2}
\|{\bar{u}}_{t}\|_{BMO(\R\times
  (0,T)}+\|{\bar{u}}_{xx}\|_{BMO(\R\times (0,T))}\leq
c \|\tilde{f}\|_{BMO(\R\times\R_{+})}.
\end{equation}
The $BMO$ theory for parabolic equations, particularly estimate
(\ref{ch3:maison-1}) is rather classical. This is due to the fact
that the solution of (\ref{maison1}) can be expressed in terms of
the heat kernel $\Gamma$ defined by:
$$\Gamma(x,t)=\left\{
\begin{aligned}
&  (4\pi \eps t)^{-1/2}e^{-\frac{x^{2}}{4\eps t}}, \quad
&\mbox{for}& \quad
t>0\\
& 0 \quad &\mbox{for}& \quad t\leq 0,
\end{aligned}
\right.$$ in the following way:
$$\bar{u}(x,t)= \int_{\R\times \R^{+}}
\Gamma(x-\xi,t-s)\tilde{f}(\xi,s)\,d\xi\,ds.$$ As a matter of fact,
it is shown in \cite{FabRivi} that $\Gamma_{xx}$ is a parabolic
Calderon-Zygmund kernel (here we are working in nonhomogeneous
metric spaces in which the variable $t$ accounts for twice the
variable $x$).  Therefore $\Gamma_{xx}:BMO\rightarrow BMO$ is a
bounded linear operator. This result is quite technical and can be
adapted from its elliptic version (see \cite[Theorem 3.4]{BraBra}).
It is less difficult to show that $\Gamma_{xx}:L^{\infty}\rightarrow
BMO$, a bounded linear operator (see for instance \cite[Lemma
3.3]{HeLaVe}).

Having (\ref{maison2}) in hands, it remains to show that
$$
 \|\tilde{f}\|_{BMO(\R\times\R_{+})}\leq
 c\left(\|f\|_{BMO(\R\times (0,T))}+m_{2I\times (0,T)}(|f|)\right),
$$
with $c>0$ independent of $T$. This can be divided into two steps:\\

\noindent {\bf Step 1. }\textsf{(Treatment of small parabolic
  cubes)}\\

\noindent We consider parabolic cubes $Q_{r}=Q_{r}(x_{0},t_{0})$,
$(x_{0},t_{0})\in \R\times \R_{+}$, with $r\leq \sqrt{T}$. Let us
estimate the  term
$$\frac{1}{|Q_{r}|}\int_{Q_{r}}|\tilde{f}-m_{Q_{r}}\tilde{f}|.$$
Assume, without loss of generality, that $T\leq t_{0}<2T$. In fact,
any other case can be done in a similar way because of the time
symmetry of the function $\tilde{f}$. Two cases can be considered.
If $r^{2}< t_{0}-T$ then the cube $Q_{r}$ lies in the strip
$\R\times (T,2T)$ and in this case
$$\frac{1}{|Q_{r}|}\int_{Q_{r}}|\tilde{f}-m_{Q_{r}}\tilde{f}|\leq
\|f\|_{BMO(\R\times(0,T))}.$$ The other case is when $r^{2}\geq
t_{0}-T$. In this case, define $Q^{a}_{r}$ and $Q^{b}_{r}$, the
above and the below parabolic cubes, as follows:
$$Q^{a}_{r}=Q_{r}(x_{0},T+r^{2})\quad \mbox{and}\quad
Q^{b}_{r}=Q_{r}(x_{0},T).$$ Since
$$T-r^{2}< t_{0}-r^{2}\leq T < t_{0} \leq T+r^{2},$$
then $Q_{r}\subset (Q^{a}_{r}\cup Q^{b}_{r})$. Moreover, we have
$|Q_{r}|=|Q^{a}_{r}|=|Q^{b}_{r}|$. We compute:
\begin{eqnarray*}
\frac{1}{|Q_{r}|}\int_{Q_{r}}|\tilde{f}-m_{Q_{r}}\tilde{f}|&\leq&
\frac{2}{|Q_{r}|}\int_{Q_{r}}
|\tilde{f}-2m_{Q^{b}_{r}}\tilde{f}+m_{Q^{a}_{r}}\tilde{f}|\\
&\leq& \frac{4}{|Q_{r}|}\int_{Q^{a}_{r}}
|\tilde{f}-m_{Q^{b}_{r}}\tilde{f}|+
\frac{4}{|Q_{r}|}\int_{Q^{b}_{r}}
|\tilde{f}-m_{Q^{b}_{r}}\tilde{f}|\\
& &
+\frac{2}{|Q_{r}|}\int_{Q^{a}_{r}}|\tilde{f}-m_{Q^{a}_{r}}\tilde{f}|+
\frac{2}{|Q_{r}|}\int_{Q^{b}_{r}}|\tilde{f}-m_{Q^{a}_{r}}\tilde{f}|.
\end{eqnarray*}
We remark (from the symmetry-in-time of the function $\tilde{f}$)
that $m_{Q^{a}_{r}}\tilde{f}=m_{Q^{b}_{r}}f$, and
$$\int_{Q^{a}_{r}}|\tilde{f}-c|=\int_{Q^{b}_{r}}|f-c|,\quad \forall
c\in \R.$$ Therefore the above inequalities give:
$$
\frac{1}{|Q_{r}|}\int_{Q_{r}}|\tilde{f}-m_{Q_{r}}\tilde{f}|\leq 16
\|f\|_{BMO(\R\times (0,T))}.
$$

\noindent {\bf Step 2. }\textsf{(Treatment of big parabolic
  cubes)}\\

\noindent Consider now parabolic cubes $Q_{r}\subset \R\times
\R_{+}$, $r>\sqrt{T}$. Suppose first that $r>1$. Because of the
symmetry-in-time of the function $\tilde{f}$, and its spatial
periodicity, we compute:
\begin{equation*}
\frac{1}{|Q_{r}|}\int_{Q_{r}}|\tilde{f}-m_{Q_{r}}\tilde{f}|\leq
\frac{2}{|Q_{r}|}\int_{Q_{r}}|\tilde{f}| \leq
\frac{2N}{|Q_{r}|}\int_{2I\times (0,T)}|f|,
\end{equation*}
where $N$ is the minimum number of domains $D$ of the form $D=
(k,k+2)\times (nT, (n+1)T)$, $k\in \Z$ and $n\in N$, that cover
$Q_{r}$. Here
$$|Q_{r}|\sim N \times |2I\times (0,T)|,\quad N>1.$$
Therefore, the above inequalities give:
$$
\frac{1}{|Q_{r}|}\int_{Q_{r}}|\tilde{f}-m_{Q_{r}}\tilde{f}|\leq c\,
m_{2I\times (0,T)}(|f|).
$$
Now suppose that $\sqrt{T}<r\leq 1$. In this case we use the fact
that $0<T_{1}\leq T$, we compute:
\begin{equation*}
\frac{1}{|Q_{r}|}\int_{Q_{r}}|\tilde{f}-m_{Q_{r}}\tilde{f}|\leq
\frac{2}{|Q_{r}|}\int_{Q_{r}}|\tilde{f}|\leq
\frac{2N}{|Q_{r}|}\int_{2I\times (0,T)}|f| \leq
\frac{N}{T^{3/2}_{1}}\int_{2I\times (0,T)}|f|.
\end{equation*}
Here $N \lesssim \frac{1}{T}$, and hence
$$
\frac{1}{|Q_{r}|}\int_{Q_{r}}|\tilde{f}-m_{Q_{r}}\tilde{f}|\leq
c(T_{1})\, m_{2I\times (0,T)}(|f|).
$$
Steps 1 and 2 give the required result. $\hfill\Box$\\

\noindent {\bf B2. Proof of Lemma \ref{3n2bnkzBTA}.} We divide the
proof
into two steps.\\

\noindent {\bf Step 1.} \textsf{(Treatment of small parabolic
  cubes)}\\

\noindent Let us consider parabolic cubes
$Q=Q_{r}(x_{0},t_{0})\subset \R\times (0,T)$ with $0<r\leq
\frac{1}{2}$. Assume, without loss of generality, that $1<x_{0}<2$
(the other cases can be treated similarly). Define the left and the
right neighbor cubes of $Q_{r}(x_{0},t_{0})$ by
$Q^{-}=Q^{-}_{r}(1-r,t_{0})$, and $Q^{+}=Q^{+}_{r}(1+r,t_{0})$
respectively. Since $2r\leq 1$, then
$$Q^{-}\subset (0,1)\times (0,T)\quad \mbox{and} \quad Q^{+}\subset
(1,2)\times (0,T).$$
Using the fact that for any function $g\in L^{1}(\Omega)$:
$$\int_{\Omega}|g-m_{\Omega}(g)|\leq 2 \int_{\Omega} |g-c|,\quad \forall
c\in \R,$$
We compute:
\begin{eqnarray}\label{Ft5Bd9856}
\frac{1}{|Q|} \int_{Q}|f^{sym}-m_{Q}(f^{sym})|&\leq & \frac{2}{|Q|}
\int_{Q} |f^{sym}+m_{Q^{+}}(f^{asym})|\nonumber\\
&\leq & \frac{2}{|Q^{-}|}
\int_{Q^{-}} |f^{sym}+m_{Q^{+}}(f^{asym})|\nonumber\\
&+& \frac{2}{|Q^{+}|}
\int_{Q^{+}} |f^{sym}+m_{Q^{+}}(f^{asym})|.
\end{eqnarray}
We know that from the properties of $f^{sym}$ and $f^{asym}$ that
$m_{Q^{+}}(f^{asym})=-m_{Q^{-}}(f^{sym})$,  and
$$f^{sym}= -f^{asym}\quad \mbox{on}\quad Q^{+},\quad \mbox{and} \quad
f^{sym}= f^{asym}\quad \mbox{on}\quad Q^{-}.$$
Using the above two inequalities in (\ref{Ft5Bd9856}), we get:
\begin{equation*}
\frac{1}{|Q|} \int_{Q}|f^{sym}-m_{Q}(f^{sym})|\leq 4
\|f^{asym}\|_{BMO(\R\times (0,T))}.
\end{equation*}

\noindent {\bf Step 2.} \textsf{(Treatment of big parabolic
  cubes)}\\

\noindent Consider parabolic cubes $Q=Q_{r}\subset \R\times (0,T)$ such
that $r>\frac{1}{2}$. In this case, we compute:
\begin{equation*}
\frac{1}{|Q|} \int_{Q}|f^{sym}-m_{Q}(f^{sym})|\leq \frac{2}{|Q|}
\int_{Q}|f^{sym}| \leq \frac{2N}{|Q|} \int_{2I\times
(0,T)}|f^{sym}|,
\end{equation*}
with
$$|Q|\sim N \times |2I\times (0,T)|,$$
therefore
$$\frac{1}{|Q|} \int_{Q}|f^{sym}-m_{Q}(f^{sym})|\leq c\,m_{2I\times
  (0,T)}(|f^{sym}|),$$
  where $c$ is a universal constant. Steps 1 and 2 directly implies the result. $\hfill{\Box}$\\

\end{section}

\noindent \textbf{Acknowledgments}. This work was supported by the
contract ANR MICA (2006-2009). The authors would like to thank
J\'{e}r\^{o}me Droniou for his valuable remarks while reading the
manuscript of the paper.

\bibliographystyle{siam}
\bibliography{biblio}

\def\cprime{$'$}
\begin{thebibliography}{10}

\bibitem{adams}
{\sc R.~A. Adams}, {\em Sobolev spaces}, Academic Press, New York-London, 1975.
\newblock Pure and Applied Mathematics, Vol. 65.

\bibitem{BraBra}
{\sc M.~Bramanti and L.~Brandolini}, {\em Estimates of {BMO} type for singular
  integrals on spaces of homogeneous type and applications to hypoelliptic
  {PDE}s}, Rev. Mat. Iberoamericana, 21 (2005), pp.~511--556.

\bibitem{BreGal}
{\sc H.~Br{\'e}zis and T.~Gallou{\"e}t}, {\em Nonlinear {S}chr\"odinger
  evolution equations}, Nonlinear Anal., 4 (1980), pp.~677--681.

\bibitem{BreWai}
{\sc H.~Br{\'e}zis and S.~Wainger}, {\em A note on limiting cases of {S}obolev
  embeddings and convolution inequalities}, Comm. Partial Differential
  Equations, 5 (1980), pp.~773--789.

\bibitem{vazquez_carrillo03}
{\sc J.~A. Carrillo and J.~L. V{\'a}zquez}, {\em Fine asymptotics for fast
  diffusion equations}, Comm. Partial Differential Equations, 28 (2003),
  pp.~1023--1056.

\bibitem{ChanKaper}
{\sc C.~Y. Chan and H.~G. Kaper}, {\em Quenching for semilinear singular
  parabolic problems}, SIAM J. Math. Anal., 20 (1989), pp.~558--566.

\bibitem{vazquez_cha02}
{\sc E.~Chasseigne and J.~L. V{\'a}zquez}, {\em Theory of extended solutions
  for fast-diffusion equations in optimal class of data. {R}adiation from
  singularities}, Arch. Ration. Mech. Anal., 164 (2002), pp.~133--187.

\bibitem{DiBook}
{\sc E.~DiBenedetto}, {\em Degenerate parabolic equations}, Universitext,
  Springer-Verlag, New York, 1993.

\bibitem{Engler}
{\sc H.~Engler}, {\em An alternative proof of the {B}rezis-{W}ainger
  inequality}, Comm. Partial Differential Equations, 14 (1989), pp.~541--544.

\bibitem{FabRivi}
{\sc E.~B. Fabes and N.~M. Rivi{\`e}re}, {\em Singular intervals with mixed
  homogeneity}, Studia Math., 27 (1966), pp.~19--38.

\bibitem{afriedman}
{\sc A.~Friedman}, {\em Partial differential equations of parabolic type},
  Prentice-Hall Inc., Englewood Cliffs, N.J., 1964.

\bibitem{GCZ}
{\sc I.~Groma, F.~F. Czikor, and M.~Zaiser}, {\em Spatial correlations and
  higher-order gradient terms in a continuum description of dislocation
  dynamics}, Acta Mater, 51 (2003), pp.~1271--1281.

\bibitem{SouGuo}
{\sc J.-S. Guo and P.~Souplet}, {\em Fast rate of formation of dead-core for
  the heat equation with strong absorption and applications to fast blow-up},
  Math. Ann., 331 (2005), pp.~651--667.

\bibitem{HayWol}
{\sc N.~Hayashi and W.~von Wahl}, {\em On the global strong solutions of
  coupled {K}lein-{G}ordon-{S}chr\"odinger equations}, J. Math. Soc. Japan, 39
  (1987), pp.~489--497.

\bibitem{HeLaVe}
{\sc M.~A. Herrero, A.~A. Lacey, and J.~J.~L. Vel{\'a}zquez}, {\em Global
  existence for reaction-diffusion systems modelling ignition}, Arch. Rational
  Mech. Anal., 142 (1998), pp.~219--251.

\bibitem{HL}
{\sc J.~R. Hirth and L.~Lothe}, {\em Theory of dislocations}, Second edition,
  Kreiger publishing company, Florida 32950, 1982.

\bibitem{IJM_PI}
{\sc H.~Ibrahim, M.~Jazar, and R.~Monneau}, {\em Dynamics of dislocation
  densities in a bounded channel. {P}art {I}: smooth solutions to a singular
  coupled parabolic system}, preprint hal-00281487, 65 pages.

\bibitem{IJM_PII}
\leavevmode\vrule height 2pt depth -1.6pt width 23pt, {\em Dynamics of
  dislocation densities in a bounded channel. {P}art {II}: existence of weak
  solutions to a singular {H}amilton-{J}acobi/parabolic strongly coupled
  system}, preprint hal-00281859, 33 pages.

\bibitem{Ib_Mo08}
{\sc H.~Ibrahim and R.~Monneau}, {\em A parabolic version of the
  {K}ozono-{T}aniuchi inequality}, in preparation.

\bibitem{KT}
{\sc H.~Kozono and Y.~Taniuchi}, {\em Limiting case of the {S}obolev inequality
  in {BMO}, with application to the {E}uler equations}, Comm. Math. Phys., 214
  (2000), pp.~191--200.

\bibitem{LSU}
{\sc O.~A. Lady{\v{z}}enskaja, V.~A. Solonnikov, and N.~N. Ural{\cprime}ceva},
  {\em Linear and quasilinear equations of parabolic type}, Translated from the
  Russian by S. Smith. Translations of Mathematical Monographs, Vol. 23,
  American Mathematical Society, Providence, R.I., 1967.

\bibitem{lieberman}
{\sc G.~M. Lieberman}, {\em Second order parabolic differential equations},
  World Scientific Publishing Co. Inc., River Edge, NJ, 1996.

\bibitem{Antonio}
{\sc A.~Maugeri}, {\em A boundary value problem for a class of singular
  parabolic equations}, Boll. Un. Mat. Ital. B (5), 17 (1980), pp.~325--339.

\bibitem{MerleZaag}
{\sc F.~Merle and H.~Zaag}, {\em Optimal estimates for blowup rate and behavior
  for nonlinear heat equations}, Comm. Pure Appl. Math., 51 (1998),
  pp.~139--196.

\bibitem{nab}
{\sc F.~R.~N. Nabarro}, {\em Theory of crystal dislocations}, Oxford, Clarendon
  Press, 1969.

\end{thebibliography}

\end{document}